\documentclass[11pt]{article}
\usepackage{amssymb,latexsym}
\usepackage{epsfig}
\usepackage{eufrak}
\usepackage{amsmath}
\usepackage{mathrsfs}
\usepackage{color}

\newtheorem{theorem}{Theorem}
\newtheorem{lemma}{Lemma}

\newtheorem{definition}{Definition}

\newcommand{\be}{\begin{equation}}
\newcommand{\ee}{\end{equation}}
\newcommand{\bee}{\begin{eqnarray*}}
\newcommand{\eee}{\end{eqnarray*}}
\newcommand{\bel}{\begin{eqnarray}}
\newcommand{\eel}{\end{eqnarray}}
\newcommand{\bec}{\begin{cases}}
\newcommand{\eec}{\end{cases}}
\newcommand{\bem}{\begin{bmatrix}}
\newcommand{\eem}{\end{bmatrix}}

\newcommand{\la}{\label}
\newcommand{\li}{\left}
\newcommand{\ri}{\right}

\newcommand{\ovl}{\overline}
\newcommand{\udl}{\underline}

\newcommand{\lc}{\lceil}
\newcommand{\rc}{\rceil}
\newcommand{\lf}{\lfloor}
\newcommand{\rf}{\rfloor}

\newcommand{\ep}{\epsilon}
\newcommand{\vep}{\varepsilon}
\newcommand{\lm}{\lambda}

\newcommand{\Up}{\Upsilon}

\newcommand{\si}{\sigma}

\newcommand{\de}{\delta}

\newcommand{\vDe}{\varDelta}

\newcommand{\ga}{\gamma}
\newcommand{\Ga}{\Gamma}
\newcommand{\vse}{\vartheta}
\newcommand{\se}{\theta}

\newcommand{\ze}{\zeta}
\newcommand{\al}{\alpha}
\newcommand{\ba}{\beta}

\newcommand{\ro}{\rho}
\newcommand{\ka}{\kappa}

\newcommand{\f}{\frac}
\newcommand{\sq}{\sqrt}
\newcommand{\cd}{\cdots}

\newcommand{\qu}{\quad}
\newcommand{\qqu}{\qquad}
\newcommand{\fa}{\forall}

\newcommand{\mscr}{\mathscr}
\newcommand{\mcal}{\mathcal}
\newcommand{\mbf}{\mathbf}

\newcommand{\wh}{\widehat}

\newcommand{\mrm}{\mathrm}
\newcommand{\bs}{\boldsymbol}

\newcommand{\ap}{\approx}

\newcommand{\LRA}{\Longleftrightarrow}
\newcommand{\sh}{\slash}

\newcommand{\tx}{\text}

\newcommand{\iy}{\infty}

\newcommand{\bed}{\begin{description}}
\newcommand{\eed}{\end{description}}
\newcommand{\bei}{\begin{itemize}}
\newcommand{\eei}{\end{itemize}}
\newcommand{\ben}{\begin{enumerate}}
\newcommand{\een}{\end{enumerate}}
\newcommand{\bib}{\bibitem}
\newcommand{\beL}{\begin{lemma}}
\newcommand{\eeL}{\end{lemma}}
\newcommand{\beT}{\begin{theorem}}
\newcommand{\eeT}{\end{theorem}}
\newcommand{\sect}{\section}

\newcommand{\bpf}{\begin{pf}}
\newcommand{\epf}{\end{pf}}
\newcommand{\bsk}{\bigskip}

\newcommand{\pfbox}{\hfill\mbox{$\Box$}}

\newenvironment{pf}{\paragraph*{Proof{\rm.}}}{\pfbox\bigskip}

\begin{document}

\title{{\bf Multistage Hypothesis Tests for the Mean of a Normal Distribution}
\thanks{The author had been previously working with Louisiana State University at Baton Rouge, LA 70803, USA,
and is now with Department of Electrical Engineering, Southern
University and A\&M College, Baton Rouge, LA 70813, USA; Email:
chenxinjia@gmail.com}}

\author{Xinjia Chen}

\date{October 2008}

\maketitle

\begin{abstract}

In this paper, we have developed new multistage tests which
guarantee prescribed level of power and are more efficient than
previous tests in terms of average sampling number and the number of
sampling operations. Without truncation, the maximum sampling
numbers of our testing plans are absolutely bounded.  Based on
geometrical arguments, we have derived extremely tight bounds for
the operating characteristic function.  To reduce the computational
complexity for the relevant integrals, we propose adaptive scanning
algorithms which are not only useful for present hypothesis testing
problem but also for other problem areas.

\end{abstract}

\sect{Introduction}

Consider a Gaussian random variable $X$ with mean $\mu$ and variance
$\si^2$.  In many applications, it is an important problem to
determine whether the mean $\mu$ is less or greater than a
prescribed value $\ga$ based on i.i.d. random samples $X_1, X_2,
\cd$ of $X$. Such problem can be put into the setting of testing
hypothesis $\mscr{H}_0 : \mu \leq \mu_0$ versus $\mscr{H}_1 : \mu >
\mu_1$ with $\mu_0 = \ga - \vep \si$ and $\mu_1 = \ga + \vep \si$,
where $\vep$ is a positive number specifying the width of the
indifference zone $(\mu_0, \mu_1)$.  It is usually required that the
size of the Type I error is no greater than $\al \in (0,1)$ and the
size of the Type II error is no greater than $\ba \in (0,1)$.  That
is,  \be \la{re1No}
 \Pr \li \{ \tx{Reject} \;
\mscr{H}_0 \mid \mu \ri \} \leq \al, \qu \fa \mu \in (-\iy, \mu_0]
\ee
 \be
 \la{re2No}
 \Pr \li \{ \tx{Accept} \; \mscr{H}_0 \mid \mu \ri
\} \leq \ba, \qu \fa \mu \in [\mu_1, \iy). \; \ee The hypothesis
testing problem described above has been extensively studied in the
framework of sequential probability ratio test (SPRT), which was
established by Wald \cite{Wald} during the period of second world
war of last century.  The SPRT
 suffers from several drawbacks.  First, the sampling number of SPRT
 is a random number which is not bounded. However, to be useful, the maximum sampling number of
 any testing plan should be bounded by a deterministic number.  Although this can be fixed
 by forced termination (see, e.g., \cite{Ghosh} and the references therein), the prescribed level of power may not be
 ensured as a result of truncation.  Second, the number of sampling
 operations of SPRT is as large as the number of samples.  In practice, it
 is usually much more economical to take a batch of samples at a
 time instead of one by one.  Third, the efficiency of SPRT is
 optimal only for the endpoints of the indifference zone. For other parametric values, the SPRT
  can be extremely inefficient.  Needless to say, a truncated version of SPRT may suffer from the same
 problem due to the partial use of the boundary of SPRT.  Third,
 when the variance $\si^2$ is not available, a weighting function
 needs to be constructed so that the testing problem can be fit into
 the framework of SPRT.  The construction of such weighting function
 is a difficult task and severely limit the efficiency of the
 resultant test plan.

 In this paper, to overcome the limitations of existing tests for the mean of a normal distribution,
  we have established a new class testing plans having the following features:
 i) The testing has a finite number of stages and thus the cost of
 sampling operations is reduced as compared to SPRT. ii) The sampling number is absolutely
 bounded without truncation. iii) The prescribed level of
 power is rigorously guaranteed.  iv) The testing is not only efficient for
 the endpoints of indifference zone, but also efficient for other parametric
 values. v) Even the variance $\si^2$ is unknown, our test plans do
 not require any weighting function.

In general, our testing plans consist of $s$ stages. For $\ell = 1,
\cd, s$, the sample size of the $\ell$-th stage is $n_\ell$. For the
$\ell$-th stage, a decision variable $\bs{D}_\ell$ is defined by
using samples $X_1, \cd, X_{n_\ell}$ such that $\bs{D}_\ell$ assumes
only three possible values $0, \; 1$ and $2$ with the following
notion:

(i) Sampling is continued until $\bs{D}_\ell \neq 0$ for some $\ell
\in \{1, \cd, s\}$. Since the sampling must be terminated at or
before the $s$-th stage, it is required that $\bs{D}_s \neq 0$.  For
simplicity of notations, we also define $\bs{D}_0 = 0$.

(ii) The null hypothesis $\mscr{H}_0$ is accepted at the $\ell$-th
stage if $\bs{D}_\ell = 1$ and $\bs{D}_i = 0$ for $1 \leq i < \ell$.

(iii) The null hypothesis $\mscr{H}_0$ is rejected at the $\ell$-th
stage if $\bs{D}_\ell = 2$ and $\bs{D}_i = 0$ for $1 \leq i < \ell$.

As will be seen in the our specific testing plans, the sample sizes
$n_1 < n_2 < \cd, n_s$ and decision variables $\bs{D}_1, \cd,
\bs{D}_s$ depend on the parameters $\al, \; \ba, \; \mu_0, \; \mu_1$
and other parameters such as the  {\it risk tuning parameter} $\ze$
and the {\it sample size incremental factor} $\ro$.  The
requirements of power can be satisfied by determining an appropriate
value of  $\ze$ via bisection search.  For this purpose, we have
derived, by a geometrical approach, readily computable bounds for
the evaluation of the operating characteristic (OC) function.

The remainder of the paper is organized as follows.  In Section 2,
we present our approach for testing the mean of a normal
distribution in the context of knowing the variance $\si^2$. In
Section 3, we describe our method for for testing the mean of a
normal distribution for situations that the variance $\si^2$ is not
available.  Section 4 discusses the evaluation of OC functions.  In
Section, we propose adaptive scanning algorithms for integration,
summation, zero finding and optimization. These new methods are
useful for our current problem and other problem areas.  Section 6
is the conclusion. All proofs of theorems are given in Appendices.

Throughout this paper, we shall use the following notations.  The
ceiling function is denoted  by $\lc . \rc$ (i.e., $\lc x \rc$
represents the smallest integer no less than $x$).  The gamma
function is denoted by $\Ga(.)$.  The inverse cosine function taking
values on $[0, \pi]$ is denoted by $\arccos (.)$. The inverse
tangent function taking values on $\li [ -\f{\pi}{2}, \f{\pi}{2} \ri
]$ is denoted by $\arctan (.)$.  We use the notation $\Pr \{ . \mid
\se \}$ to indicate that the associated random samples $X_1, X_2,
\cd$ are parameterized by $\se$. The parameter $\se$ in $\Pr \{ .
\mid \se \}$  may be dropped whenever this can be done without
introducing confusion. The other notations will be made clear as we
proceed.

\section{ Testing the Mean of a Normal Distribution with Known
Variance}

For $\de \in (0, 1)$, let $\mcal{Z}_{\de} > 0 $ be the critical
value of a normal distribution with zero mean and unit variance,
i.e., $\Phi(\mcal{Z}_{\de}) = \f{1}{\sq{2 \pi}} \int_{
\mcal{Z}_{\de} }^\iy e^{- \f{x^2}{2} } d x = \de$.  In situations
that the variance $\si^2$ is known, our testing plan, developed in
\cite{Chen_TH},  is described as follows.

\beT \la{THMT9} Let $\ze > 0$ and $\ro > 0$. Let $n_1 < n_2 < \cd <
n_s$ be the ascending arrangement of all distinct elements of
{\small $\li \{ \li \lc \f{ (\mcal{Z}_{\ze \al} + \mcal{Z}_{\ze
\ba})^2 }{ 4 \vep^2 } (1 + \ro)^{i - \tau} \ri \rc: i = 1, \cd, \tau
\ri \}$}, where $\tau$ is a positive integer.  Define $a_\ell = \vep
\sq{n_\ell} - \mcal{Z}_{\ze \ba}, \; b_\ell = \mcal{Z}_{\ze \al} -
\vep \sq{n_\ell}$ for $\ell = 1, \cd, s - 1$, and {\small $a_s = b_s
= \f{ \mcal{Z}_{\ze \al} - \mcal{Z}_{\ze \ba} }{2}$}.  Define{\small
\[
 \ovl{X}_{n_\ell} = \f{ \sum_{i = 1}^{n_\ell} X_i } {
n_\ell }, \qqu T_\ell = \f{ \sq{n_\ell} \; (\ovl{X}_{n_\ell} - \ga)
} {\si}, \qqu \bs{D}_\ell = \bec 1 &
\tx{for} \; T_\ell \leq a_\ell,\\
2 & \tx{for} \; T_\ell > b_\ell,\\
 0 & \tx{else}  \eec
\]}
for $\ell = 1, \cd, s$.  Then, both (\ref{re1No}) and (\ref{re2No})
are guaranteed  provided that $0 < \ze \leq \f{1}{\tau}$. Moreover,
the OC function $\Pr \li \{ \tx{Accept} \; \mscr{H}_0 \mid \mu \ri
\}$ is monotonically decreasing with respect to $\mu \in (-\iy,
\mu_0) \cup (\mu_1, \iy)$.

\eeT

To compute tight bounds for the OC function, we have the following
result.

\beT \la{OC_bounds_known} Let $U$ and $V$ be independent Gaussian
random variables with zero mean and variance unity.  Define {\small
\bee \varphi(\se, \ze, \al, \ba) & = &  \Phi \li ( \sq{n_1} \se -
b_1 \ri ) + \sum_{\ell = 2}^{s} \Pr \li \{ b_\ell - \sq{n_\ell} \se
\leq U
\leq k_\ell V - \sq{n_\ell} \se +  \sq{\f{n_\ell}{n_{\ell - 1}}} b_{\ell - 1} \ri \}\\
&   & -  \; \sum_{\ell = 2}^{s}  \Pr \li \{b_\ell - \sq{n_\ell} \se
\leq U \leq k_\ell V - \sq{n_\ell} \se +  \sq{\f{n_\ell}{n_{\ell -
1}}} a_{\ell - 1}  \ri \}
 \eee}
with $k_\ell = \sq{ \f{n_\ell }{n_{\ell - 1} } - 1 }, \; \ell = 2,
\cd, s$.  Then, $\Pr \{ \tx{Accept} \;  \mscr{H}_0 \mid \mu = \se
\si + \ga \} > 1 - \varphi(\se, \ze, \al, \ba)$ for any $\se \in
(-\iy , - \vep]$ and $\Pr \{ \tx{Accept} \;  \mscr{H}_0 \mid \mu =
\se \si + \ga \} < \varphi(-\se, \ze, \ba, \al)$ for any $\se \in
[\vep, \iy)$. \eeT

\bsk See Appendix A for a proof.

As can be seen from the proof of Theorem 2, we have {\small
$\sum_{\ell = 1}^s \Pr \{ \bs{D}_{\ell - 1} = 0, \; \bs{D}_\ell = 2
\mid \mu_0 \} = \varphi (\f{\mu_0 - \ga}{\si}, \ze, \al, \ba)$} and
{\small $\sum_{\ell = 1}^s \Pr \{ \bs{D}_{\ell - 1} = 0, \;
\bs{D}_\ell = 1 \mid \mu_1 \} = \varphi(\f{\ga - \mu_1}{\si}, \ze,
\ba, \al)$}. By making use of such results and a bisection search
method, we can determine an appropriate value of $\ze$ so that both
(\ref{re1No}) and (\ref{re2No}) are guaranteed.

With regard to the distribution of sample number $\mbf{n}$, we have,
for $\ell =  1, \cd, s - 1$, \bee \Pr \{ \mathbf{n} > n_\ell \} &
\leq  &  \Pr \li \{ a_\ell  < T_\ell  \leq  b_\ell \ri \} = \Pr \li
\{ T_\ell \leq  b_\ell \ri \} -  \Pr
 \li \{ T_\ell \leq  a_\ell \ri \}\\
&  = &  \Pr \li \{ U + \se \sq{n_\ell} \leq b_\ell  \ri \} - \Pr \li
\{ U + \se \sq{n_\ell} \leq a_\ell \ri \} =  \Phi \li ( b_\ell  -
\sq{n_\ell} \se \ri ) - \Phi \li ( a_\ell -  \sq{n_\ell} \se \ri ),
\eee where $U$ is a Gaussian random variable with zero mean and unit
variance.

\section{Testing the Mean of a Normal Distribution with Unknown
Variance}

For $\de \in (0, 1)$, let $t_{n, \de}$ be the critical value of
Student's $t$-distribution with $n$ degrees of freedom. Namely,
$t_{n, \de}$ is a number satisfying
\[
\int_{t_{n, \de}}^\iy \f{ \Ga ( \f{n + 1}{2} ) } { \sq{n \pi} \; \Ga
( \f{n}{2}) }  \li ( 1 + \f{x^2}{n} \ri )^{ - \f{ n + 1 } {2} } =
\de.
\]
In situations that the variance $\si^2$ is unknown, our testing
plan, developed in \cite{Chen_TH}, is described as follows.

\beT \la{THMT10} Let $\ze > 0$ and $\ro > 0$. Let $n^*$ be the
minimum integer $n$ such that {\small $t_{n - 1, \ze \al}  + t_{n -
1, \ze \ba} \leq 2 \vep \sq{n - 1}$}. Let $n_1 < n_2 < \cd < n_s$ be
the ascending arrangement of all distinct elements of {\small $ \{
\lc n^* \; (1 + \ro)^{i - \tau} \rc: i = 1, \cd, \tau  \}$}, where
$\tau$ is a positive integer.  Define $a_\ell = \vep \sq{n_\ell - 1}
- t_{n_\ell - 1, \ze \ba}, \; b_\ell = t_{n_\ell - 1, \ze \al} -
\vep \sq{n_\ell - 1}$ for $\ell = 1, \cd, s - 1$, and {\small $a_s =
b_s = \f{ t_{n_s - 1, \ze \al} - t_{n_s - 1, \ze \ba} }{2}$}. Define
{\small
\[ \ovl{X}_{n_\ell} = \f{ \sum_{i = 1}^{n_\ell} X_i } { n_\ell },
\qu \wh{\si}_{n_\ell} = \sq{ \f{\sum_{i=1}^{n_\ell} (X_i -
\ovl{X}_{n_\ell})^2}{n_\ell - 1} }, \qu  \wh{T}_\ell = \f{
\sq{n_\ell} (\ovl{X}_{n_\ell} - \ga)}{\wh{\si}_{n_\ell}}, \qu
\bs{D}_\ell = \bec 1 &
\tx{for} \; \wh{T}_\ell \leq a_\ell,\\
2 & \tx{for} \; \wh{T}_\ell > b_\ell,\\
 0 & \tx{else}  \eec
\]}
for $\ell = 1, \cd, s$.  Then, both (\ref{re1No}) and (\ref{re2No})
are guaranteed if $\ze > 0$ is sufficiently small. Moreover, the OC
function $\Pr \li \{ \tx{Accept} \; \mscr{H}_0 \mid \mu \ri \}$ is
monotonically decreasing with respect to $\mu \in (- \iy, \mu_0)
\cup (\mu_1, \iy)$.

\eeT

To obtain tight bounds for the OC function, the following result is
useful.

\beT \la{OC_T} Let $U, \; V$ and $Y_\ell, \; Z_\ell, \; \ell = 2,
\cd, s$ be independent random variables such that $U, \; V$ possess
identical normal distributions with zero mean and unit variance and
that $Y_\ell, \; Z_\ell$ possess chi-square distributions of
$n_{\ell - 1} - 1$ and $n_{\ell} - n_{\ell - 1} - 1$ degrees of
freedom respectively.  Define {\small $k_\ell = \sq{
\f{n_\ell}{n_{\ell -1}} - 1  }, \; c_\ell = \f{a_\ell } { \sq{n_\ell
- 1} }$} and {\small $d_\ell = \f{b_\ell  } { \sq{n_\ell - 1} }$}
for $\ell = 1, \cd, s$. Define $\mcal{P}(\se, \ze, \al, \ba) =
\sum_{\ell = 1}^s \mcal{P}_\ell$ where {\small $\mcal{P}_1 = \Pr \li
\{ \wh{T}_1 > b_1 \ri \}$} and {\small \[ \mcal{P}_\ell = \bec \Pr
\li \{ d_\ell \sq{ V^2 + Y_\ell + Z_\ell } <  U  + \sq{n_\ell} \se
\leq k_\ell V + d_{\ell -
1} \sq{ \f{n_\ell Y_\ell} {n_{\ell - 1}} } \ri \}\\
 - \; \Pr \li \{ d_\ell \sq{
V^2 + Y_\ell + Z_\ell } <  U  + \sq{n_\ell} \se \leq k_\ell V  +
c_{\ell -
1} \sq{ \f{n_\ell Y_\ell} {n_{\ell - 1}} } \ri \} & \tx{for} \; d_\ell \geq 0,\\

\Pr \li \{ a_{\ell - 1} < \wh{T}_{\ell - 1}  \leq  b_{\ell - 1} \ri
\} +  \Pr \li \{ | d_\ell | \sq{ V^2 + Y_\ell + Z_\ell } \leq U -
\sq{n_\ell} \se < k_\ell V - d_{\ell
- 1} \sq{ \f{n_\ell Y_\ell} {n_{\ell - 1}} } \ri \}\\
 -  \; \Pr \li \{ | d_\ell | \sq{ V^2 + Y_\ell + Z_\ell } \leq U - \sq{n_\ell} \se < k_\ell V - c_{\ell
- 1} \sq{ \f{n_\ell Y_\ell} {n_{\ell - 1}} } \ri \} & \tx{for} \;
d_\ell < 0 \eec
\]}
for $\ell = 2, \cd, s$.  Then, $\Pr \{ \tx{Accept} \; \mscr{H}_0
\mid \mu = \se \si + \ga \} \geq 1 - \mcal{P} (\se, \ze, \al, \ba)$
for any $\se \in (-\iy , - \vep]$ and $\Pr \{ \tx{Accept} \;
\mscr{H}_0 \mid \mu = \se \si + \ga \} \leq \mcal{P} (-\se, \ze,
\ba, \al)$ for any $\se \in [\vep, \iy)$. \eeT

\bsk

See Appendix B for a proof. As can be seen from the proof of Theorem
4,  we have $\sum_{\ell = 1}^s \Pr \{ \bs{D}_{\ell - 1} = 0, \;
\bs{D}_\ell = 2 \mid \mu_0 \} = \mcal{P} (\f{\mu_0 - \ga}{\si}, \ze,
\al, \ba)$ and $\sum_{\ell = 1}^s \Pr \{ \bs{D}_{\ell - 1} = 0, \;
\bs{D}_\ell = 1 \mid \mu_1 \} = \mcal{P} (\f{\ga - \mu_1}{\si}, \ze,
\ba, \al)$. By making use of such results and a bisection search
method, we can determine an appropriate value of $\ze$ so that both
(\ref{re1No}) and (\ref{re2No}) are guaranteed.

With regard to the distribution of the sample number $\mbf{n}$, we
have {\small $\Pr \{ \mathbf{n}
> n_\ell \} < \Pr \{ a_\ell < \wh{T}_\ell
 \leq b_\ell \}$} for $\ell = 1, \cd, s-1$, where the probability
 can be expressed in terms of the well-known non-central $t$-distribution.

\section{Evaluation of OC Functions}

In this section, we shall demonstrate that the evaluation of OC
functions of tests described in preceding discussion can be reduced
to the computation of the probability of a certain domain including
two independent standard Gaussian variables.  In this regard, our
first general result is as follows.

\beT Let $U$ and $V$ be two independent Gaussian random variables
with zero mean and unit variance. Let $\mscr{D}$ be a
two-dimensional convex domain which contains the origin $(0,0)$.
Suppose the set of boundary points of $\mscr{D}$ can be expressed as
$\mscr{B} = \{ (r, \phi ): r = \mcal{B} (\phi), \; \phi \in \mscr{A}
\}$ in polar coordinates, where $\mcal{B} (\phi)$ is a Riemann
integrable function on set $\mscr{A}$. Then,
\[
\Pr \{ (U, V) \in \mscr{D} \} = 1 - \f{1}{2 \pi} \int_{\mscr{A}}
\exp \li ( - \f{ \mcal{B}^2 (\phi) } { 2 } \ri ) d \phi.
\]
\eeT

See Appendix C for a proof.  For situations that the domain does not
contain the origin $(0,0)$, we need to introduce the concept of {\it
visibility} for boundary points of a two-dimensional domain
$\mscr{D}$.  The intuitive notion of such concept is that a boundary
point of $\mscr{D}$ is visible if it can be seen by an observer at
the origin. The precise definition is as follows.

\begin{definition}
A  boundary point, $(u, v)$, of domain $\mscr{D}$ is said to be
visible if $\{ (q u, q v) : 0 < q < 1 \} \cap \mscr{D}$ is empty.
Otherwise, such a boundary point is said to be invisible.
\end{definition}

Based on the concept of visibility, we have derived the following
general result.

\beT Let $U$ and $V$ be two independent Gaussian random variables
with zero mean and unit variance. Let $\mscr{D}$ be a
two-dimensional convex domain which does not contain the origin
$(0,0)$. Suppose the set of visible boundary points of $\mscr{D}$
can be expressed as $\mscr{B}_{\mrm{v}} = \{ (r, \phi ): r =
\mcal{B}_{\mrm{v}} (\phi), \; \phi \in \mscr{A}_{\mrm{v}} \}$ in
polar coordinates, where $\mcal{B}_{\mrm{v}} (\phi)$ is a Riemann
integrable function on set $\mscr{A}_{\mrm{v}}$.  Suppose the set of
invisible boundary points of $\mscr{D}$ can be expressed as
$\mscr{B}_{\mrm{i}} = \{ (r, \phi ): r = \mcal{B}_{\mrm{i}} (\phi),
\; \phi \in \mscr{A}_{\mrm{i}} \}$ in polar coordinates, where
$\mcal{B}_{\mrm{i}} (\phi)$ is a Riemann integrable function on set
$\mscr{A}_{\mrm{i}}$.  Then,
\[
\Pr \{ (U, V) \in \mscr{D} \} = \f{1}{2 \pi} \li [
\int_{\mscr{A}_{\mrm{v}}} \exp \li ( - \f{ \mcal{B}_{\mrm{v}}^2
(\phi) } { 2 } \ri ) d \phi - \int_{\mscr{A}_{\mrm{i}}} \exp \li ( -
\f{ \mcal{B}_{\mrm{i}}^2 (\phi) } { 2 } \ri ) d \phi \ri ].
\]
\eeT

\bsk

See Appendix D for a proof. As can be seen from Theorem 2, the
evaluation of OC functions of test plans designed for the case that
the variance $\si^2$ is known can be reduced to the computation of
probabilities of the form $\Pr \{ h \leq U \leq k V + g \}$.  For
fast computation of such probabilities, we have derived, based on
Theorems 5 and 6, the following result. \beT \la{Int_Thm} Let $k >
0$. Let $U$ and $V$ be independent Gaussian random variables with
zero mean and unit variance. Define {\small $\Psi_h (\phi) = \f{1}{2
\pi} \exp \li ( - \f{ h^2 } { 2 \cos^2 \phi } \ri ), \; \Psi_{g,k}
(\phi) = \f{1}{2 \pi} \exp \li ( - \f{ g^2 } { 2 (1 + k^2) \cos^2
\phi } \ri ), \; \phi_k = \arctan \li ( k \ri )$} and {\small
$\phi_R = \arctan \li ( \f{h - g}{k h} \ri )$}. Then, \[ \Pr \{ h
\leq U \leq k V + g \} = \bec \int^{\pi + \phi_k + \phi_R }_{\pi \sh
2} \Psi_{g,k} (\phi) \; d \phi -
\int^{\pi + \phi_R }_{\pi \sh 2} \Psi_h (\phi) \; d \phi  & \tx{for  $\max (g,h) < 0$},\\

1 - \int^{\pi + \phi_R}_{\pi \sh 2}  \Psi_h (\phi) \; d \phi
- \int^{3 \pi \sh 2}_{\phi_k + \phi_R}  \Psi_{g,k} (\phi) \; d \phi & \tx{for  $h \leq 0 \leq g$},\\

\int^{\pi \sh 2}_{\phi_R} \Psi_h (\phi) \; d \phi - \int^{\pi \sh
2}_{\phi_k + \phi_R} \Psi_{g,k} (\phi) \; d \phi  & \tx{else}.
  \eec
\]
\eeT

\bsk

See Appendix E for a proof. As can be seen from Theorem 4, the
evaluation of OC functions of test plans designed for the case that
the variance $\si^2$ is unknown can be reduced to the computation of
probabilities of the type {\small $\Pr \li \{ \lm \sq{ V^2 + Y + Z }
\leq U - \vse < k V + \varpi \sq{ Y } \ri \}$} with $\lm > 0$, where
$Y$ and $Z$ are chi-square random variables independent with $U$ and
$V$.  The evaluation of such probabilities is described as follows.

Define multivariate functions $\ovl{P} (\udl{y}, \ovl{z})$ and
$\udl{P} (\udl{y}, \ovl{z})$ so that {\small
\[ \ovl{P} (\udl{y}, \ovl{z}) = \bec \Pr \li \{ \lm \sq{V^2 +
\udl{y} + \udl{z}} \leq U - \vse \leq k V + \varpi \sq{ \ovl{y} } \ri \} & \tx{if} \; \varpi \geq 0,\\
\Pr \li \{ \lm \sq{V^2 + \udl{y} + \udl{z}} \leq U - \vse \leq k V +
\varpi \sq{ \udl{y} } \ri \} & \tx{if} \; \varpi < 0 \eec
\]
\[
\udl{P} (\udl{y}, \ovl{z}) = \bec \Pr \li \{ \lm \sq{ V^2 + \ovl{y}
+ \ovl{z}} \leq U - \vse \leq  k V + \varpi \sq{
\udl{y} } \ri \} & \tx{if} \; \varpi \geq 0,\\
\Pr \li \{ \lm \sq{V^2 + \ovl{y} + \ovl{z}} \leq U - \vse \leq k V +
\varpi \sq{ \ovl{y} }  \ri \} & \tx{if} \; \varpi < 0 \eec
\]}
for $0 < \udl{y} \leq \ovl{y}, \; 0 < \udl{z} \leq \ovl{z}$.  Then,
$\Pr \{ \lm \sq{ V^2 + Y + Z} \leq U - \vse \leq  k V + \varpi
\sq{Y}, \; Y \in [\udl{y}, \ovl{y}] , \; Z \in [ \udl{z},  \ovl{z} ]
\}$ is smaller than $\Pr \li \{ Y \in [ \udl{y}, \ovl{y} ] \ri \}
\times \Pr \li \{ Z \in [ \udl{z}, \ovl{z} ] \ri \} \times \ovl{P}
(\udl{y},  \ovl{z})$ and is greater than $\Pr \li \{ Y \in [
\udl{y}, \ovl{y} ] \ri \} \times \Pr \li \{ Z \in [ \udl{z}, \ovl{z}
] \ri \} \times \udl{P} (\udl{y}, \ovl{z})$.  For any $\ep \in (0,
1)$, we can determine, via bisection search, positive numbers
$y_{\mrm{min}} < y_{\mrm{max}}$ and $z_{\mrm{min}} < z_{\mrm{max}}$
such that $\Pr \{ Y < y_{\mrm{min}} \} < \f{\ep}{4}, \; \Pr \{ Y
> y_{\mrm{max}} \} < \f{\ep}{4}, \; \Pr \{ Z < z_{\mrm{min}} \} < \f{\ep}{4}$ and $\Pr \{ Z > z_{\mrm{max}} \}
< \f{\ep}{4}$. By partitioning the set $\{ (y,z) : y \in [
y_{\mrm{min}},  y_{\mrm{max}}], \; z \in [ z_{\mrm{min}},
z_{\mrm{max}} ] \}$ as sub-domains $\{ (y,z) : y \in [ \udl{y}_i,
\ovl{y}_i ], \; z \in [  \udl{z}_i,  \ovl{z}_i ] \}, \; i = 1, \cd,
m$ and evaluating $\ovl{P}_i = \Pr \{ Y \in [ \udl{y}_i, \ovl{y}_i ]
\} \times \Pr \li \{ Z \in [ \udl{z}_i, \ovl{z}_i ] \ri \} \times
\ovl{P} (\udl{y}_i, \ovl{z}_i)$ and $\udl{P}_i = \Pr \{ \udl{y}_i
\leq Y \leq \ovl{y}_i \} \times \Pr \li \{ Z \in [ \udl{z}_i,
\ovl{z}_i ] \ri \} \times \udl{P} (\udl{y}_i, \ovl{z}_i)$ for $i =
1, \cd, m$, we have
 {\small \[ \sum_{i} \udl{P}_i < \Pr \li \{ \lm \sq{ V^2 + Y + Z}
\leq U - \vse \leq  k V + \varpi \sq{ Y } \ri \} < \ep + \sum_{i}
\ovl{P}_i. \] The bounds can be refined by further partitioning the
sub-domains. For efficiency, we can split the sub-domain with the
largest gap between the upper bound $\ovl{P}_i$ and lower bound
$\ovl{P}_i$ in every additional partition.  It can be seen that the
probabilities like $\ovl{P} (\udl{y}_i, \ovl{z}_i)$ and $\udl{P}
(\udl{y}_i, \ovl{z}_i)$ are of the same type as $\Pr \{ (U, V) \in
\mscr{D} \}$,  where {\small $\mscr{D} =  \{ (u,v) : \sq{\lm v^2 +
h} \leq u  - \vse \leq k v + g  \}$} with $k > 0, \; \lm > 0, \; h
\geq 0$ and $k^2 \neq \lm$.  For fast computation of such
probabilities, we have derived, based on Theorems 5 and 6, the
following results.

\beT Define {\small $\vDe = h (k^2 - \lm) + \lm g^2, \; u_A = \f{
\lm g - k \sq{ \vDe } } { \lm - k^2 } + \vse, \; u_B = \f{ \lm g + k
\sq{ \vDe }   } { \lm - k^2 } + \vse, \; v_A =  \f{ g k - \sq{ \vDe}
} { \lm - k^2 }, \; v_B = \f{ g k + \sq{ \vDe} } { \lm - k^2 }, \;
\phi_A = \arccos \li (  \f{u_A}{\sq{u_A^2 + v_A^2} } \ri ), \;
\phi_B = \arccos \li ( \f{u_B}{\sq{u_B^2 + v_B^2}} \ri ), \;
\phi_{\mrm{m}} = \arctan \li ( \sq{ \f{h}{ \lm |\vse^2 - h|  } } \ri
), \; \phi_{\lm}  = \arctan \li ( \f{1}{\sq{\lm} } \ri )$}, {\small
$\phi_k = \arctan (k), \; \Psi_{\vse, g, k} ( \phi)  = \f{1}{2 \pi}
\exp \li ( - \f{(\vse + g)^2} {2 (1 + k^2) \cos^2 \phi } \ri )$} and
{\small
\[
\Upsilon_{\vse, \lm, h} ( \phi) = \f{1}{2 \pi} \exp \li ( - \f{
(\vse^2 - h)^2 } { 2 \li [\vse \cos \phi + \sq{ (h + \lm h - \lm
\vse^2 ) \cos^2 \phi + \lm (\vse^2 - h) } \ri ]^2} \ri ).
\]}
Then, {\small \[
 \Pr \{ (U, V)
\in \mscr{D} \} = \bec I_{\mrm{np}} & \tx{for} \; k^2 < \lm, \; g > \sq{h}, \; \vDe \geq 0,\\
I_{\mrm{pp}} & \tx{for} \; k^2 < \lm, \; 0 < g \leq \sq{h}, \; \vDe \geq 0,\\
I_{\mrm{n}} & \tx{for} \; k^2 > \lm, \;g k > \sq{ \vDe},\\
I_{\mrm{p}} & \tx{for} \; k^2 > \lm, \;g k \leq \sq{ \vDe},\\
0 & \tx{else}  \eec
\]}
where {\small \[
I_{\mrm{np}} = \bec I_{\mrm{np},1} & \tx{for} \; \vse + \f{h}{u_B - \vse} \geq 0,\\
I_{\mrm{np},2} & \tx{for} \; \vse + \f{h}{u_B - \vse} < 0 \leq \vse + \f{h}{u_A - \vse},\\
I_{\mrm{np},3} & \tx{for} \; \vse + \f{h}{u_A - \vse} < 0 \leq \vse + \sq{h},\\
I_{\mrm{np},4} & \tx{for} \; \vse + \sq{h} < 0 \leq \vse + g,\\
I_{\mrm{np},5} & \tx{for} \; \vse + g < 0 \eec \qqu \qu
I_{\mrm{n}} = \bec I_{\mrm{n},1} & \tx{for} \; \vse \geq 0,\\
I_{\mrm{n},2} & \tx{for} \; \vse < 0 \leq \vse + \f{h}{u_A - \vse},\\
I_{\mrm{n},3} & \tx{for} \; \vse + \f{h}{u_A - \vse} < 0 \leq \vse + \sq{h},\\
I_{\mrm{n},4} & \tx{for} \; \vse + \sq{h} < 0 \leq \vse + g, \\
I_{\mrm{n},5} & \tx{for} \; \vse + g < 0 \eec
\]}
{\small \[
I_{\mrm{pp}} = \bec I_{\mrm{pp},1} & \tx{for} \; \vse + \f{h}{u_B - \vse} \geq 0,\\
I_{\mrm{pp},2} & \tx{for} \; \vse + \f{h}{u_B - \vse} < 0 \leq \vse + \f{h}{u_A - \vse},\\
I_{\mrm{pp},3} & \tx{for} \; \vse + \f{h}{u_A - \vse} < 0 \eec \qqu
\qu
I_{\mrm{p}} = \bec I_{\mrm{p},1} & \tx{for} \; \vse \geq 0,\\
I_{\mrm{p},2} & \tx{for} \; \vse < 0 \leq \vse + \f{h}{u_A - \vse},\\
I_{\mrm{p},3} & \tx{for} \; \vse + \f{h}{u_A - \vse} < 0 \eec \qqu
\qu
\]}
with {\small \bee I_{\mrm{np},1} & = & \int_{\pi-\phi_A}^{\pi +
\phi_B} \Upsilon (\phi) d \phi -  \int_{\phi_k
-\phi_A }^{ \phi_k + \phi_B } \Psi (\phi) d \phi,\\
I_{\mrm{np},2} & = & \int^{\pi + \phi_{\mrm{m}} }_{\pi - \phi_A}
\Upsilon (\phi) d \phi - \int^{\phi_{\mrm{m}} }_{\phi_B} \Upsilon
(\phi) d \phi - \int_{\phi_k -\phi_A }^{\phi_k + \phi_B } \Psi
(\phi) d \phi,\\
I_{\mrm{np},3} & = &  \int_{\pi - \phi_{\mrm{m}} }^{\pi +
\phi_{\mrm{m}} } \Upsilon (\phi) d \phi - \int^{\phi_{\mrm{m}}
}_{\phi_B} \Upsilon (\phi) d \phi
 - \int_{\phi_A}^{\phi_{\mrm{m}} } \Upsilon (\phi) d \phi  -  \int_{\phi_k
-\phi_A }^{\phi_k + \phi_B} \Psi (\phi) d \phi,\\
I_{\mrm{np},4} & = &  1 - \int_{\phi_k -
\phi_A }^{ \phi_k + \phi_B } \Psi (\phi)  d \phi - \int_{\phi_B}^{2 \pi - \phi_A} \Upsilon (\phi) d \phi,\\
I_{\mrm{np},5} & = & \int_{\phi_k + \phi_B  }^{\phi_k - \phi_A + 2
\pi} \Psi (\phi) d \phi - \int_{\phi_B}^{2 \pi - \phi_A} \Upsilon
(\phi) d \phi, \eee} {\small \bee
 I_{\mrm{n},1} &  = &  \int_{\pi - \phi_A}^{\pi +
\phi_{\lm} } \Upsilon (\phi) d \phi - \int^{ \f{\pi}{2}
}_{\phi_k - \phi_A } \Psi (\phi) d \phi ,\\
I_{\mrm{n},2} & = &  \int_{\pi - \phi_A}^{\pi + \phi_{\mrm{m}} }
\Upsilon (\phi) d \phi - \int_{\phi_{\lm} }^{\phi_{\mrm{m}} }
\Upsilon (\phi) d \phi -  \int^{ \f{\pi}{2} }_{\phi_k - \phi_A }
\Psi (\phi) d \phi,\\
I_{\mrm{n},3} & = &  \int_{\pi - \phi_{\mrm{m}} }^{\pi +
\phi_{\mrm{m}} } \Upsilon (\phi) d \phi - \int_{\phi_{\lm}
}^{\phi_{\mrm{m}} } \Upsilon (\phi) d \phi  -
\int_{\phi_A}^{\phi_{\mrm{m}} } \Upsilon (\phi) d \phi -
\int^{ \f{\pi}{2} }_{\phi_k - \phi_A} \Psi (\phi) d \phi,\\
I_{\mrm{n},4} & = & 1 - \int^{ \f{\pi}{2} }_{\phi_k - \phi_A} \Psi
(\phi) d \phi -
\int_{\phi_{\lm} }^{2 \pi - \phi_A} \Upsilon (\phi) d \phi,\\
I_{\mrm{n},5} & = & \int_{ \f{\pi}{2} }^{\phi_k - \phi_A + 2 \pi }
\Psi (\phi) d \phi - \int_{\phi_{\lm} }^{2 \pi - \phi_A} \Upsilon
(\phi) d \phi, \\
 I_{\mrm{pp},1} & = & \int_{\pi + \phi_A}^{\pi +
\phi_B} \Upsilon (\phi) d \phi - \int_{\phi_k
+ \phi_A}^{\phi_k + \phi_B} \Psi (\phi) d \phi,\\
 I_{\mrm{pp},2} & = &
\int_{\pi + \phi_A}^{\pi + \phi_{\mrm{m}} } \Upsilon (\phi) d \phi -
\int_{\phi_B}^{\phi_{\mrm{m}} } \Upsilon (\phi) d \phi -
\int_{\phi_k + \phi_A}^{\phi_k + \phi_B} \Psi (\phi) d \phi,\\
I_{\mrm{pp},3} & = & \int_{\phi_k + \phi_B}^{\phi_k + \phi_A} \Psi
(\phi) d \phi - \int_{\phi_B}^{\phi_A} \Upsilon (\phi) d \phi,\\
 I_{\mrm{p},1} &  = & \int_{ \f{\pi}{2} }^{\phi_k + \phi_A } \Psi
(\phi) d \phi + \int_{\pi + \phi_A}^{\pi + \phi_{\lm} } \Upsilon
(\phi)
d \phi,\\
I_{\mrm{p},2} & = & \int_{ \f{\pi}{2} }^{\phi_k + \phi_A } \Psi
(\phi)  d \phi + \int_{\pi + \phi_A}^{\pi + \phi_{\mrm{m}} }
\Upsilon (\phi) d \phi
- \int_{\phi_{\lm} }^{\phi_{\mrm{m}} } \Upsilon (\phi) d \phi,\\
I_{\mrm{p},3} & = & \int_{ \f{\pi}{2} }^{\phi_k + \phi_A } \Psi
(\phi) d \phi - \int_{\phi_{\lm} }^{\phi_A} \Upsilon (\phi) d \phi.
 \eee}
\eeT

\bsk

See Appendix F for a proof. In Theorem 8, for simplicity of
notations, we have abbreviated $\Psi_{\vse, g, k} ( \phi)$ and
$\Upsilon_{\vse, \lm, h} ( \phi)$ as $\Psi ( \phi)$ and $\Upsilon (
\phi)$ respectively.

\section{Adaptive Scanning Algorithms}

As can be seen from last section, we need to frequently evaluate
integrals involving functions like $\Up(.)$ and $\Psi(.)$. Clearly,
there are no closed-form solutions for this type of integrals.
Although existing numerical integration method can be applied to
obtain approximations for such integrals, the accuracy of
integration is not clearly known.  Since our concern is the risk of
making wrong decisions in hypothesis testing, the quantification of
integration is crucial.  Motivated by this consideration, we have
developed an adaptive scanning method for fast integration.
Moreover, we have extended the method to summation, zero finding and
optimization.

\subsection{Integration of Continuous Functions}

The integrals involved in hypothesis testing can be addressed in the
general framework of computing $I (a, b) = \int_a^b f (x) dx$ by a
numerical method.  Existing methods are quadrature rules.

A quadrature rule is an approximation of the definite integral of a
function, usually stated as a weighted sum of a function values at
specified points within the domain of integrations. More formally, a
quadrature rule proceeds as follows:

(i) Partition the interval $[a, b]$ by grid points $a = x_0 < x_1 <
\cd < x_n = b$.

(ii) Evaluate $f(x_i), \; i = 0, 1, \cd, n$.

(iii) Construct an estimate $\wh{I} (a, b)$ for $I (a, b)$ as a
weighted sum of $f(x_i), i = 0, 1, \cd, n$.

Well known quadrature rules are rectangle rule, trapezium rule,
Simpson's rule, Romberg's method, Gaussian quadrature rule,
Clenshaw-Curtis quadrature rule, Newton-Cote formula, Richardson
extrapolation, etc.

One of the most frequently used method is the composite Simpson's
rule.  Suppose that the interval $[a, b]$ is split up in $n$
subintervals, with $n$ an even number. Then, the composite Simpson's
rule is given by
\[
\int_a^b f(x) dx \ap \f{h}{3} \li [ f (x_0) + 2 \sum_{j =
1}^{\f{n}{2} - 1} f (x_{2 j}) + 4 \sum_{j = 1}^{\f{n}{2}} f (x_{2 j
- 1}) + f (x_{n}) \ri ],
\]
where $x_j = a + j h, \; j = 0, 1, \cd, n - 1, n$ and $h = \f{b -
a}{n}$; in particular, $x_0 = a$ and $x_n = b$.

It is widely recognized that an assessment of the accuracy is an
essential part of any numerical method.  Specifically, given $\vep
> 0$, a crucial question is how to ensure
\[
| \wh{I} (a, b) - I (a, b) | \leq \vep?
\]
The error committed by the composite Simpson's rule is bounded (in
absolute value) by
\[
\f{h^4}{180} (b - a) \max_{\ze \in [a, b]} \li | f^4 (\ze) \ri |.
\]
In Simpson's rule, it is not clear how to choose the step length. If
the step length is too small, the computation is too slow. On the
other hand, a large step length may cause intolerable error of the
computation.  Although the error bound can be expressed in terms of
the fourth derivative, to guarantee the accuracy, we need to bound
the fourth derivative over the whole integration range $[a, b]$. The
bounding is not easy and can be extremely conservative.

It is not hard to see that other quadrature rules suffer similar
drawbacks as the Simpson's rule.   To overcome such drawbacks,  we
propose a new approach so that the accuracy requirement can be
rigorously guaranteed under mild conditions. A salient feature of
our approach is that, instead of partition the interval $[a, b]$, we
sequentially and adaptively perform integration over subintervals of
the overall interval.  For each subinterval, making use of
derivative information, we force the integration to meet a certain
accuracy requirement. Starting from the left endpoint of interval
$[a, b]$, we determine an initial $[u_1, v_1]$ with $u_1 = a$ such
that the difference between $I (u_1,v_1) = \int_{u_1}^{v_1} f(x) dx$
and its estimate $\wh{I} (u_1, v_1)$ is no greater than $\f{\vep}{b
- a}(v_1 - u_2)$.  Then, we determine next subinterval $[u_2, v_2]$
as the form
\[ u_2 = v_1, \qqu v_2 = \min \{b, v_1 + (v_1 - u_1) 2^{\ell} \},
\]
with $\ell$ taken as the largest integer no greater than $1$ to
ensure that the difference between $I (u_2,v_2) = \int_{u_2}^{v_2}
f(x) dx$ and its estimate $\wh{I} (u_2, v_2)$ is no greater than
$\f{\vep}{b - a}(v_2 - u_2)$.  For $i > 1$, given interval $[u_i,
v_i]$,  we determine next subinterval $[u_{i+1}, v_{i+1}]$ as the
form
\[ u_{i+1} = v_i, \qqu v_{i+1} = \min \{b, v_i + (v_i - u_i) 2^{\ell} \},
\]
with $\ell$ taken as the largest integer no greater than $1$ to
ensure that the difference between $I (u_{i+1},v_{i+1}) =
\int_{u_{i+1}}^{v_{i+1}} f(x) dx$ and its estimate $\wh{I} (u_{i+1},
v_{i+1})$ is no greater than $\f{\vep}{b - a}(v_{i+1} - u_{i+1})$.
We repeat this process until $v_i = b$ for some $i$.  Finally, the
overall estimate $\wh{I} (a, b)$ for $I(a, b)$ is given by
\[
\wh{I} (a, b) = \sum_{i} \wh{I} (u_i, v_i),
\]
which ensures that
\[
| \wh{I} (a, b)  - I (a, b) | \leq \sum_{i} \li | \wh{I} (u_i, v_i)
- I (u_i, v_i) \ri | \leq \sum_{i}  \f{\vep}{ b - a} (v_i - u_i) =
\vep.
\]
Since the above process of integration is like scanning the interval
of integration, we call the method as Adaptive Scanning Algorithm
(ASA).  The adaptive nature of the algorithm can be seen from the
dynamic choice of the length of subinterval $[u_i, v_i]$.   To
formally describe our ASA, let $I (u, v) = \int_{u}^{v} f(x) dx$ and
$\wh{I}(u, v)$ be an estimate of $I(u,v)$ for $a \leq u \leq v \leq
b$.  Assume that $| \wh{I} (u, v) - I (u, v) | \to 0$ as $| u - v|
\to 0$.  Let $\eta = \f{\vep}{b - a}$. Assume that we have a method
for testing the truth of $| \wh{I} (u, v) - I (u, v) | \leq \eta (v
- u)$ without knowing $I (u, v)$.  Our ASA proceeds as follows.

\bsk

\begin{tabular} {|l |}
\hline

$\diamond$  Choose initial step length $\vDe$ as a positive number less than $\f{b - a}{2}$. \\

$\diamond$ Let $\wh{I} (a, b)
\leftarrow 0, \; \eta \leftarrow \f{\vep}{b - a}$ and $u \leftarrow a$.\\

$\diamond$  While $u + \vDe < b$, do the following:\\

       $\indent$ $\star$ Let $\tx{st} \leftarrow 0$ and $\ell \leftarrow
       2$;\\

        $\indent$  $\star$ While $\tx{st} = 0$, do the following:\\

               $\indent$ $\indent$ $*$ Let $\ell \leftarrow \ell - 1$ and $\vDe \leftarrow 2^\ell \vDe$.\\
               $\indent$ $\indent$ $*$ If $u + \vDe < b$, let $v \leftarrow u + \vDe$.  Otherwise, let $v \leftarrow b$.\\

               $\indent$ $\indent$ $*$ Evaluate $\wh{I} (u, v)$. \\

               $\indent$ $\indent$ $*$ Test the truth of $| \wh{I} (u, v) - I (u, v) | \leq \eta (v -
               u)$\\

               $\indent$ $\indent$  $\;\;$ without knowledge of $I (u,
               v)$.\\

               $\indent$ $\indent$ $*$ If $| \wh{I} (u, v) - I (u, v) | \leq \eta (v - u)$ is true, \\

               $\indent$ $\indent$ $\; \;$ let $\wh{I} (a, b)
\leftarrow \wh{I} (a, b) + \wh{I} (u, v)$ and $\tx{st} \leftarrow 1, \; u \leftarrow   v$.\\

$\diamond$  Return $\wh{I} (a, b)$ as an estimate for $I(a, b)$.\\
\hline
\end{tabular}

\bsk

Under the assumption that $| \wh{I} (u, v) - I (u, v) | \to 0$ as $|
u - v| \to 0$, it can be readily shown that $| \wh{I} (a, b) - I (a,
b) | \leq \vep$ is guaranteed after execution of the algorithm. This
is because $| \wh{I} (a, b) - I (a, b) |$ is no greater than the
summation of $| \wh{I} (u, v) - I (u, v) |$ over all subintervals
$(u,v)$ generated to cover $[a, b]$.

As can be seen from the description of ASA, a critical issue is to
construct $\wh{I} (u, v)$ and test the truth of $| \wh{I} (u, v) - I
(u, v) | \leq \eta (v - u)$ without any knowledge of $I (u, v)$. Our
general method for addressing this issue is as follows.  Let
$\udl{I} (u, v)$ and $\ovl{I} (u, v)$ be lower and upper bounds of
$I(u, v)$ respectively.  Namely,  $\udl{I} (u, v) \leq I (u, v) =
\int_{u}^{v} f(x) dx \leq \ovl{I} (u, v) $. Assume that $ \ovl{I}
(u, v) - \udl{I} (u, v) \to 0$ as $| u - v| \to 0$. In many cases,
the lower and upper bounds can be obtained from Taylor series
expansion formula. To test the truth of $| \wh{I} (u, v) - I (u, v)
| \leq \eta (v - u)$, we propose to make use of the following
relationship \be \la{imp}
 \ovl{I} (u, v) - \eta (v - u) \leq \wh{I} (u, v) \leq
\udl{I} (u, v) + \eta (v - u) \qu \Longrightarrow \qu | \wh{I} (u,
v) - I (u, v) | \leq \eta (v - u). \ee  To construct estimate
$\wh{I} (u, v)$ for $I(u, v)$, we recommend to take $\wh{I} (u, v) =
\f{1}{2} [ \udl{I} (u, v) + \ovl{I} (u, v) ]$ or
\[
\wh{I} (u, v) =  \f{v - u}{6} \li [ f(u) + 4 f  ( \f{u + v}{2} ) +
f(v) \ri ]
\]
based on Simpson's approximation rule.

Assuming that the first derivative $f^\prime (x)$ of $f(x)$ exists
for all $x \in [a, b]$, making use of (\ref{imp}), Taylor's series
expansion formula, and Simpson's approximation rule, we have derived
the following methods in Theorem \ref{rule} for constructing $\wh{I}
(u, v)$ and testing the truth of $| \wh{I} (u, v) - I (u, v) | \leq
\eta (v - u)$ without any knowledge of $I (u, v)$.

\beT \la{rule} Let $u, v$ and $w$ be three real numbers such that $w
- u = v - w = h > 0$. Let $I(u, v) = \int_{u}^{v} f(x) dx$ and
$\wh{I} (u, v) = \f{h}{3} [ f(u) + 4 f(w) + f(v) ]$. Then, the
following statements hold true.

(I) $|\wh{I} (u, v) - I(u,v)| \leq \eta (v - u)$ provided that
\[
3 \ovl{\ka}   - \f{6 \eta}{h} \leq \f{f(u) + f(v) - 2 f(w)}{ h} \leq
3 \udl{\ka} + \f{ 6 \eta } {h},
\]
where $\udl{\ka} = \f{1}{2} [ \min_{x \in [u, w]} f^\prime (x) +
\min_{x \in [w, v]} f^\prime (x)]$ and $\ovl{\ka} = \f{1}{2} [
\max_{x \in [u, w]} f^\prime (x) + \max_{x \in [w, v]} f^\prime
(x)]$.

(II) $|\wh{I} (u, v) - I(u,v)| \leq \eta (v - u)$ provided that
$f(x)$ is a concave function of $x \in [u, v]$ and that
\[
- \f{ 12 \eta }{h} \leq  \f{f(u) + f(v) - 2 f(w)}{h} \leq \f{3}{4}
\li [ f^\prime (v) - f^\prime (u) \ri ] + \f{ 12 \eta }{h}.
\]

(III) $|\wh{I} (u, v) - I(u,v)| \leq \eta (v - u)$ provided that
$f(x)$ is a convex function of $x \in [u, v]$ and that
\[
\f{3}{4} \li [ f^\prime (v) - f^\prime (u) \ri ]  - \f{ 12 \eta }{h}
\leq  \f{f(u) + f(v) - 2 f(w)}{h} \leq \f{ 12 \eta }{h}.
\]

\eeT

\bsk

In the case that the convexity or concavity of $f(x)$ are hard to
determine, one may compute the bounds of the first derivative of
$f(x)$ and apply statement (I) of Theorem \ref{rule} to ASA.  For
example, the derivatives of elliptical functions can be easily
bounded, and thus one can use statement (I) for the purpose of
integration.

The applications of statements (I) and (II) of Theorem \ref{rule}
depend on the convexity or concavity of $f(x)$.  To determine the
convexity or concavity of $f(x)$, we can find the inflexion points
from the equation $f^{\prime \prime} (x) = 0$, which can frequently
be reduced to a quadratic equation of $x$. Specially, this is true
for normal distribution, Gamma distribution, Beta distribution,
Student's $t$-distribution, and $F$-distribution, etc. Once the
inflexion points are obtained, the interval of integration can be
decomposed as subintervals so that $f(x)$ is completely convex or
concave in each subinterval.

\subsection{Summation of Discrete Functions}

In parallel to the problem of computing $I(a, b) = \int_a^b f(x)
dx$,  a similar problem is the computation of the discrete summation
$S(a, b) = \sum_{k = a}^b f(k)$, where $a, b, k$ are integers. Let
$\udl{S} (u, v)$ and $\ovl{S} (u, v)$ be the lower and upper bounds
of $S(u, v) = \sum_{k = u}^v f(k)$ respectively. We can easily
modify the ASA of integration for computing $S(a, b)$ as follows.

\bsk

\begin{tabular} {|l |}
\hline

$\diamond$  Choose initial step length $\vDe$ as a positive integer less than $\f{b - a}{2}$. \\

$\diamond$ Let $\wh{S} (a, b)
\leftarrow 0, \; \eta \leftarrow \f{\vep}{b - a + 1}$ and $u \leftarrow a$.\\

$\diamond$  While $u + \vDe < b$, do the following:\\

       $\indent$ $\star$ Let $\tx{st} \leftarrow 0$ and $\ell \leftarrow
       2$;\\

        $\indent$  $\star$ While $\tx{st} = 0$, do the following:\\

               $\indent$ $\indent$ $*$ Let $\ell \leftarrow \ell - 1$ and $\vDe \leftarrow \lc 2^\ell \vDe \rc$.\\
               $\indent$ $\indent$ $*$ If $u + \vDe < b$, let $v \leftarrow u + \vDe$.  Otherwise, let $v \leftarrow b$.\\

               $\indent$ $\indent$ $*$ If $u + 1 < v$, evaluate $\udl{S} (u, v)$ and $\ovl{S} (u, v)$. \\

               $\indent$ $\indent$ $*$ If $u + 1 < v$ and $\ovl{S} (u, v) - \udl{S} (u, v) \leq 2 \eta (v - u + 1)$, \\

               $\indent$ $\indent$ $\; \;$ let $\wh{S} (a, b)
\leftarrow \wh{S} (a, b) + \f{1}{2} [ \ovl{S} (u, v) - \udl{S} (u, v) ]$ and $\tx{st} \leftarrow 1$.\\

               $\indent$ $\indent$ $*$ If $u + 1 = v$,  let $\wh{S} (a, b)
\leftarrow \wh{S} (a, b) + f(u) + f(v)$ and $\tx{st} \leftarrow 1$.\\

               $\indent$ $\indent$ $*$ If $\tx{st} = 1$, let $u \leftarrow   v +
               1$.\\

$\diamond$  Return $\wh{S} (a, b)$ as an estimate for $S(a, b)$.\\
\hline
\end{tabular}

\bsk

Clearly, $| \wh{S} (a, b) - S (a, b) | \leq \vep$ is guaranteed
after the execution of the above algorithm.  A key routine is to
calculate the lower and upper bounds of $S(u, v) = \sum_{k = u}^v
f(k)$.  For this purpose, we have established in \cite{Chen_EST} the
following results.

\beT  Let $u < v$ be two integers. Define {\small $r_u = \f{ f(u +
1) }{ f(u) }, \; r_v = \f{ f(v - 1) } { f (v)}, \; r_{u, v} = \f{
f(u) }{ f(v) }$} and {\small $j = u + \f{ v - u - (1 -  r_{u,v}) (1
- r_v )^{-1} } { 1 + r_{u, v} ( 1 - r_u ) ( 1 - r_v )^{-1} }$}.
Define {\small $\al (i) = (i + 1 - u) \li [ 1 + \f{(i - u) (  r_u  -
1 ) }{2}  \ri ]$} and {\small $\ba (i) = (v - i) \li [ 1 + \f{( v -
i - 1) ( r_v - 1 ) }{2}  \ri ]$}.  The following statements hold
true:

 (I): If $f(k+1) - f(k) \leq f(k) - f(k-1)$  for $u < k < v$, then
\be \la{concavity} \f{ (v - u + 1) [f(u) + f(v)] }{2}  \leq
\sum_{k=u}^v f(k) \leq \alpha (i) f(u) + \beta (i)  f(v) \ee for $u
< i < v$. The minimum gap between the lower and upper bounds is
achieved at $i$ such that $\lf j \rf \leq i \leq \lc j \rc$.

(II): If $f(k+1) - f(k) \geq f(k) - f(k-1)$  for $u < k < v$, then
\[
\f{ (v - u + 1) [f(u) + f(v)] }{2}  \geq \sum_{k=u}^v f(k) \geq
\alpha (i) f(u) + \beta (i)  f(v)
\]
for $u < i < v$. The minimum gap between the lower and upper bounds
is achieved at $i$ such that $\lf j \rf \leq i \leq \lc j \rc$.

\eeT

\bsk

To investigate conditions like $f(k+1) - f(k) \leq f(k) - f(k-1)$ or
$f(k+1) - f(k) \geq f(k) - f(k-1)$, we can find the inflexion points
from equation $f(k+1) - f(k) = f(k) - f(k-1)$, which in many cases
can be reduced to a quadratic equation of $k$. Specially, this is
true for binomial distribution, negative binomial distribution,
Poisson distribution and hyper-geometrical distribution, etc.  Once
the inflexion points are obtained, we can decompose the range of
summation as subsets so that $f(k)$ is completely convex or concave
in each subset.

\subsection{Zero Finding}

To determine the convexity or concavity of $f(x)$, we need to find
zeros of the second derivative $f^{\prime \prime} (x)$.  For a
function like $\Up (.)$, there exits no analytic solution. Motivated
by such situation,  we propose a general method for finding the
zeros of $f(x)$ for $x \in [a, b]$.  Sine the zeros can be obtained
consecutively, this problem can be reduced to the following generic
problem:

Suppose that $f(a) < 0$ and $f(x)$ is continuous for $x \in [a, b]$.
Determine whether $f(x)$ has at least one root in $[a, b]$.  In the
case that $f(x)$ has at least one root in $[a, b]$, find the
smallest root $x \in [a, b]$ such that $f(x) = 0$.

Assume that, for any interval $[u, v] \subseteq [ a, b]$, it is
possible to compute an upper bound $g (u, v)$ such that $f(x) \leq
g(u, v)$ for any $x \in [u, v]$ and that the upper bound converges
to $f(x)$ as the interval width $v - u$ tends to $0$.  Let $\eta >
0$ be an extremely small number, i.e. $\eta = 10^{-15}$. Our
algorithm for zero finding proceeds as follows:

\bsk

\begin{tabular} {|l |}
\hline

$\diamond$  Choose initial step length $\vDe$ as a number between
$\eta$ and $\f{b-a}{2}$.\\

$\diamond$  Let $F \leftarrow 0, \; T \leftarrow 0$ and $a
\leftarrow u$.\\

$\diamond$  While $F = T = 0$, do the following:\\

      $\indent$  $\star$ Let $\tx{st} \leftarrow 0$ and $\ell \leftarrow
       2$;\\

      $\indent$ $\star$ While $\tx{st} = 0$, do the following:\\

               $\indent$ $\indent$ $*$ Let $\ell \leftarrow \ell - 1$ and $\vDe \leftarrow \vDe
               2^\ell$.\\

              $\indent$ $\indent$  $*$ If $u + \vDe < b$, let $v \leftarrow u + \vDe$ and $T
               \leftarrow 0$.   Otherwise, let $v \leftarrow b$ and $T \leftarrow
               1$.\\

              $\indent$ $\indent$  $*$ If $g(u, v) < 0$, let $\tx{st} \leftarrow 1$ and $u \leftarrow
               v$.\\

                $\indent$ $\indent$ $*$ If $\vDe < \eta$, let $\tx{st} \leftarrow 1$ and $F \leftarrow
                1$.\\

$\diamond$  If $F = 1$, return $x = \f{u + v}{2}$ as the smallest
root in $[a, b]$ such that $f(x) = 0$. \\
 $\;$ Otherwise if $F = 0$,
declare that  $f(x)$ has no root on $[a,
b]$.\\

\hline
\end{tabular}

\bsk

The above algorithm declares $x = \f{u + v}{2}$ as an estimate of
the smallest root based on the observation that $f(x) < 0$ for all
$x \in [a, u]$ and that $g(u, v) \geq 0$.  Since $v - u < \eta \ap
0$ and $g(u, v) \to f (\f{u + v}{2})$ as $v - u \to 0$, it is
reasonable to believe that the smallest root is close to $\f{u +
v}{2}$. In the case that $f(x)$ has more than one roots in $[a, b]$,
the above algorithm can be repeatedly used to find all the zeros.

It should be noted that this algorithm is actually adapted from our
{\it Adaptive Maximum Checking Algorithm} (AMCA) established in
\cite{Chen_EST}.

\subsection{Finding Maximum}

Clearly, finding the zeros of function $f(x)$ is closely related to
the problem of finding the minimum or maximum of $f(x)$.   Our AMCA
can be adapted for finding the maximum of $f(x)$ for $x \in [a, b]$.

From our previous paper \cite{Chen_EST}, it can be seen that our
AMCA is a computational method to determine whether a function
$f(x)$ is smaller than a prescribed number for every value of $x$ in
interval $[a, b]$.  Suppose that we have a lower bound $L$ and an
upper bound $U$ for $\max_{x \in [a, b]} f(x)$. Then, we can apply
our AMCA and a bisection search method to determine the exact value
of $\max_{x \in [a, b]} f(x)$.  One way to find a lower bound $L$ is
to compute $n$ values of $f(x)$ and take the maximum as $L$. Once a
lower bound $L$ is obtained, one can find an upper bound $U$ as the
form $U = L 2^k$, where the positive number $k$ can be determined as
the minimum integer  by our AMCA such that $L 2^k > \max_{x \in [a,
b]} f(x)$.  Of course, there are some other methods for finding $L$
and $U$.

\sect{Conclusion}

In this paper, we have developed new multistage sampling schemes for
testing the mean of a normal distribution. Our sampling schemes have
absolutely bounded number of samples.  Our test plans are
significantly more efficient than previous tests, while rigorously
guaranteeing  prescribed level of power. In contrast to existing
tests, our test plans involve no probability ratio and weighting
function. The evaluation of operating characteristic functions of
our tests can be readily accomplished by using tight bounds derived
from a geometrical approach.  We have established adaptive scanning
methods for integration, summation, zero finding and optimization,
which are useful for our current problem and other fields.

\appendix

\section{Proof of Theorem 2}

To show Theorem 2, the following lemma is useful.

\beL \la{thmm2} Let $m < n$ be two positive integers.  Let $X_1,
X_2, \cd, X_n$ be i.i.d. normal random variables with common mean
$\mu$ and variance $\si^2$.  Let {\small $\ovl{X}_k = \f{
\sum_{i=1}^k X_i }{ k }$} for $k = 1, \cd, n$. Let {\small
$\ovl{X}_{m, n} = \f{ \sum_{i=m + 1}^n X_i }{ n - m }$}. Define
{\small \[ U = \f{\sq{n}(\ovl{X}_n - \mu)}{\si}, \;
 V = \sq{\f{m(n-m)}{n}} \f{\ovl{X}_m -  \ovl{X}_{m,n}}{\si}, \;
Y = \f{1}{\si^2} \sum_{i=1}^m \li ( X_i - \ovl{X}_m \ri )^2, \; Z =
\f{1}{\si^2} \sum_{i=m+1}^n \li ( X_i - \ovl{X}_{m,n} \ri )^2.
\]}
Then, $U, V, Y, Z$ are independent random variables such that both
$U$ and $V$ are normally distributed with zero mean and variance
$1$, $Y$ possesses a chi-square distribution of degree $m-1$, and
$Z$ possesses a chi-square distribution of degree $n-m-1$. Moreover,
$\sum_{i = 1}^n (X_i - \ovl{X}_n)^2 = \si^2 ( Y + Z + V^2)$.

\eeL

\bpf Observing that {\small $R_1 = \f{ \sq{m} (\ovl{X}_m - \mu) }{
\si}$} and {\small $R_2 = \f{ \sq{n-m} (\ovl{X}_{m,n} - \mu)
}{\si}$} are independent Gaussian random variables with zero mean
and unit variance and that $U, \; V$ can be obtained from $R_1, \;
R_2$ by an orthogonal transformation {\small \[
\bem U\\
V \eem = \bem \sq{ \f{m}{n} } &  \sq{ \f{n-m}{n} }\\
\sq{ \f{n-m}{n} } &  - \sq{ \f{m}{n} }
 \eem \bem R_1\\
R_2 \eem,
\]}
we have that $U$ and $V$ are also independent Gaussian random
variables with zero mean and unit variance. Since $R_1, \; R_2, \;
Y, \; Z$ are independent, we have that $U, V, Y, Z$ are independent.
For simplicity of notations, let $S_n = \sum_{i = 1}^n (X_i -
\ovl{X}_n)^2$ and $S_{m,n} = \sum_{i = m}^n (X_i -
\ovl{X}_{m,n})^2$. Using identity $S_n = \sum_{i=1}^n X_i^2 - n
\ovl{X}_n^2$, we have {\small $\sum_{i=1}^m X_i^2 = S_m + m
\ovl{X}_m^2, \; \sum_{i=m+1}^n X_i^2 = S_{m,n} + (n-m)
\ovl{X}_{m,n}^2$} and {\small \bee S_n & = & \sum_{i=1}^n X_i^2 - n
\ovl{X}_n^2
 =  \sum_{i=1}^m X_i^2 + \sum_{i=m+1}^n X_i^2 - n \li [ \f{m \ovl{X}_m + (n-m) \ovl{X}_{m,n}}{n} \ri ]^2\\
& = & S_m + m \ovl{X}_m^2 + S_{m,n} + (n-m) \ovl{X}_{m,n}^2 - n \li
[ \f{m \ovl{X}_m + (n-m) \ovl{X}_{m,n}}{n} \ri
]^2\\
& = & S_m + S_{m,n} + \f{ m (n-m)  } { n } ( \ovl{X}_m -
\ovl{X}_{m,n}
)^2\\
& = &  \sum_{i=1}^m (X_i - \ovl{X}_m)^2  + \sum_{i=m+1}^n (X_i -
\ovl{X}_{m,n})^2 + \f{m(n-m)}{n} ( \ovl{X}_m -  \ovl{X}_{m,n} )^2 =
\si^2 \li ( Y + Z +  V^2 \ri ). \eee}

\epf

Now we are in a position to prove the theorem.   By Lemma 1 and some
algebraic operations, we have
\[
U + \sq{ \f{n-m}{m} } V  = \f{\sq{n} (\ovl{X}_m - \mu) }{\si}, \qqu
\f{(\ovl{X}_m - \mu) }{\si} = \f{1}{\sq{n}} \li ( U + \sq{
\f{n-m}{m} } V  \ri ),
\]
\[ \f{\sq{m} (\ovl{X}_m - \ga) }{\si} = \sq{\f{m}{n}} \li ( U + \sq{n} \se + \sq{ \f{n-m}{m} } V  \ri
), \qqu \f{\sq{n}(\ovl{X}_n - \ga)}{\si} = U + \sq{n} \se.
\]
For $\ell = 1$, we have $\Pr \{ \tx{Reject} \; \mscr{H}_0, \;
\mathbf{n} = n_1 \mid \mu = \se \si + \ga \} = \Pr \{ \bs{D}_\ell =
2 \mid \mu = \se \si + \ga \}  \leq  \Pr \li \{ T_1
> b_1 \ri \} = \Pr \li \{ U + \sq{n_1} \se > b_1 \ri \} = \Phi (
\sq{n_1} \se - b_1 )$ for any $\se \in ( - \iy,  - \vep]$.  For $1 <
\ell \leq s$, since $a_{\ell - 1} \leq b_{\ell - 1}$, we have
{\small \bee \Pr \{ \tx{Reject} \; \mscr{H}_0, \; \mathbf{n} =
n_\ell \mid \mu = \se \si + \ga \}
 & < & \Pr \{ \bs{D}_{\ell -
1} = 0, \; \bs{D}_\ell = 2 \mid \mu = \se \si + \ga \}\\
&  = & \Pr \li \{ a_{\ell - 1} < T_{\ell - 1} \leq b_{\ell - 1},
 \;\; T_\ell > b_\ell \ri \}\\
& = &  \Pr \li \{ T_{\ell - 1} \leq b_{\ell - 1}, \;\; T_\ell >
b_\ell \ri \} -  \; \Pr \li \{ T_{\ell - 1} \leq a_{\ell - 1}, \;\;
T_\ell > b_\ell \ri \}\\
& = &  \Pr \li \{ \sq{\f{n_{\ell - 1}}{n_\ell}} \li ( U +
\sq{n_\ell} \se + k_\ell V \ri ) \leq b_{\ell - 1}, \;\; U +
\sq{n_\ell} \se  > b_\ell \ri \}\\
&  & -  \; \Pr \li \{ \sq{\f{n_{\ell - 1}}{n_\ell}} \li ( U +
\sq{n_\ell} \se + k_\ell V  \ri ) \leq a_{\ell - 1}, \;\; U +
\sq{n_\ell} \se > b_\ell \ri \}\\
& = & \Pr \li \{ b_\ell - \sq{n_\ell} \se \leq U
\leq k_\ell V - \sq{n_\ell} \se +  \sq{\f{n_\ell}{n_{\ell - 1}}} b_{\ell - 1} \ri \}\\
&   & -  \; \Pr \li \{b_\ell - \sq{n_\ell} \se \leq U \leq k_\ell V
- \sq{n_\ell} \se +  \sq{\f{n_\ell}{n_{\ell - 1}}} a_{\ell - 1}  \ri
\} \eee} for any $\se \in ( - \iy,  - \vep]$.  It follows that $\Pr
\{ \tx{Accept} \; \mscr{H}_0, \; \mathbf{n} = n_\ell \mid \mu = \se
\si + \ga \} = 1 - \sum_{\ell = 1}^s \Pr \{ \tx{Reject} \;
\mscr{H}_0, \; \mathbf{n} = n_\ell \mid \mu = \se \si + \ga \} > 1 -
\varphi ( \se, \ze, \al, \ba)$ for any $\se \in ( - \iy,  - \vep]$.
By symmetry, we have $\Pr \{ \tx{Accept} \; \mscr{H}_0, \;
\mathbf{n} = n_\ell \mid \mu = \se \si + \ga \} < \varphi ( - \se,
\ze, \ba, \al)$ for any $\se \in [\vep, \iy )$.  This completes the
proof of the theorem.

\section{Proof of Theorem 4}

By Lemma 1, we have
\[
\f{ \wh{T}_{\ell - 1} }{\sq{n_{\ell - 1} -1}} = \sq{ \f{n_{\ell -
1}} {n_\ell} } \f{ U + \sq{n_\ell} \se + k_\ell V  } { \sq{Y_\ell}
}, \qqu \f{ \wh{T}_\ell }{\sq{n_\ell -1}} = \f{U + \sq{n_\ell} \se}{
\sq{V^2 + Y_\ell + Z_\ell} }
\]
for $1 < \ell \leq s$.  We shall focus on the case of $\mu \leq \ga
- \vep \si$, since the case of $\mu \leq \ga + \vep \si$ can be
dealt with symmetrically. For $\ell = 1$,  we have $\Pr \{
\tx{Reject} \; \mscr{H}_0, \; \mathbf{n} = n_1 \} \leq \mcal{P}_1$
for any $\se \in ( - \iy,  - \vep]$.  For $1 < \ell \leq s$, we have
{\small \bee \Pr \{ \tx{Reject} \; \mscr{H}_0, \; \mathbf{n} =
n_\ell \mid \mu = \se \si + \ga \} & < &  \Pr \{ \bs{D}_{\ell -
1} = 0, \; \bs{D}_\ell = 2 \mid \mu = \se \si + \ga \}\\
&  =  & \Pr \li \{ a_{\ell - 1}  < \wh{T}_{\ell - 1}  \leq b_{\ell -
1}, \; \f{ \wh{T}_\ell }{\sq{n_\ell -1}}  > d_\ell \ri \}\\
& = & \Pr \li \{ a_{\ell - 1}
   < \wh{T}_{\ell - 1}  \leq  b_{\ell - 1}, \;  \f{U + \sq{n_\ell} \se}{ \sq{V^2 + Y_\ell + Z_\ell} } > d_\ell \ri
   \}
  \eee} for any $\se \in ( - \iy,  - \vep]$.  In the case of $d_\ell \geq  0$, since $c_{\ell - 1} \leq d_{\ell -
1}$, it is evident that {\small  \bee &   & \Pr \li \{ a_{\ell - 1}
   < \wh{T}_{\ell - 1}  \leq  b_{\ell - 1}, \;  \f{U + \sq{n_\ell} \se}{ \sq{V^2 + Y_\ell + Z_\ell} } > d_\ell \ri
   \}\\
 & = & \Pr \li \{ c_{\ell - 1} \sq{
\f{n_\ell Y_\ell} {n_{\ell - 1}} }
   < U + \sq{n_\ell} \se + k_\ell V  \leq  d_{\ell - 1}
   \sq{ \f{n_\ell Y_\ell} {n_{\ell - 1}} }, \;  \f{U + \sq{n_\ell} \se}{ \sq{V^2 + Y_\ell + Z_\ell} } > d_\ell \ri
   \} =  \mcal{P}_\ell
  \eee} for any $\se \in ( - \iy,  - \vep]$.  In the case of $d_\ell < 0$, we have {\small \bee  &  & \Pr \li
  \{ a_{\ell - 1} < \wh{T}_{\ell - 1}  \leq  b_{\ell - 1}, \;  \f{U + \sq{n_\ell} \se}{ \sq{V^2 + Y_\ell + Z_\ell} } > d_\ell \ri \}\\
& = & \Pr \li \{ a_{\ell - 1} < \wh{T}_{\ell - 1}  \leq  b_{\ell -
1} \ri \}   - \Pr \li \{ a_{\ell - 1} < \wh{T}_{\ell - 1}  \leq
b_{\ell - 1}, \;  \f{U + \sq{n_\ell} \se}{ \sq{V^2 + Y_\ell + Z_\ell} } \leq d_\ell \ri \}\\
& = & \Pr \li \{ a_{\ell - 1} < \wh{T}_{\ell - 1}  \leq  b_{\ell -
1} \ri \}\\
&   &   - \Pr \li \{ - d_{\ell - 1} \sq{ \f{n_\ell Y_\ell} {n_{\ell
- 1}} } < U - \sq{n_\ell} \se + k_\ell V  \leq  - c_{\ell - 1}
   \sq{ \f{n_\ell Y_\ell} {n_{\ell - 1}} }, \;  \f{U - \sq{n_\ell} \se}{ \sq{V^2 + Y_\ell + Z_\ell} } \geq - d_\ell \ri
   \} =  \mcal{P}_\ell  \eee} for any $\se \in ( - \iy,  - \vep]$.  It follows that $\Pr
\{ \tx{Accept} \; \mscr{H}_0, \; \mathbf{n} = n_\ell \mid \mu = \se
\si + \ga \} = 1 - \sum_{\ell = 1}^s \Pr \{ \tx{Reject} \;
\mscr{H}_0, \; \mathbf{n} = n_\ell \mid \mu = \se \si + \ga \} > 1 -
\mcal{P} ( \se, \ze, \al, \ba)$ for any $\se \in ( - \iy,  - \vep]$.
By symmetry, we have $\Pr \{ \tx{Accept} \; \mscr{H}_0, \;
\mathbf{n} = n_\ell \mid \mu = \se \si + \ga \} < \mcal{P} ( - \se,
\ze, \ba, \al)$ for any $\se \in [\vep, \iy )$.  This completes the
proof of the theorem.

\section{Proof of Theorem 5}

Without loss of any generality, we can assume that $\mscr{A}
\subseteq [0, 2 \pi]$ for any convex domain $\mscr{D}$ which
contains the origin $(0,0)$.  Let $\mscr{A}_* = [0, 2 \pi] \setminus
\mscr{A}$. Since $\Pr \{ (U, V) \in \mscr{D} \} = \f{1}{2 \pi} \int
\int_{ (u, v) \in \mscr{D} } \exp \li ( - \f{u^2 + v^2} {2}\ri ) du
dv$, using polar coordinates, we have \bee 2 \pi \Pr \{ (r, \phi)
\in \mscr{D} \} & = & \int_{\mscr{A}} \li [ \int_{r = 0 }^{\mcal{B}
(\phi)} \exp \li ( - \f{r^2} {2}\ri ) r dr \ri ] d \phi +
\int_{\mscr{A}_*} \li [ \int_{r = 0 }^\iy \exp \li ( - \f{r^2}
{2}\ri
) r dr \ri ] d \phi\\
& = & \int_{\mscr{A}} \li [ 1 - \exp \li ( - \f{\mcal{B}^2 (\phi)}
{2} \ri
) \ri ] d \phi + \int_{\mscr{A}_*}  d \phi\\
& = & \int_{\mscr{A} \cup \mscr{A}_*}  d \phi -  \int_{\mscr{A}}
\exp \li ( -
\f{\mcal{B}^2 (\phi)} {2} \ri ) d \phi\\
& = & 2 \pi - \int_{\mscr{A}} \exp \li ( - \f{\mcal{B}^2 (\phi)} {2}
\ri ) d \phi, \eee from which the theorem immediately follows.

\section{Proof of Theorem 6}

Without loss of any generality, we can assume that
$\mscr{A}_{\mrm{i}} \subseteq \mscr{A}_{\mrm{v}}$ for any convex
domain $\mscr{D}$ which does not contain the origin $(0,0)$.  Hence,
we can write $\mscr{D} = \mscr{D}^\prime \cup \mscr{D}^{\prime
\prime}$ with $\mscr{D}^\prime = \{ (r, \phi): \mcal{B}_{\mrm{v}}
(\phi) \leq r \leq \mcal{B}_{\mrm{i}} (\phi), \; \phi \in
\mscr{A}_{\mrm{i}} \}$ and $\mscr{D}^{\prime \prime} = \{ (r, \phi):
r \geq \mcal{B}_{\mrm{v}} (\phi), \; \phi \in \mscr{A}_{\mrm{v}}
\setminus \mscr{A}_{\mrm{i}} \}$, where $(r, \phi)$ represents polar
coordinates.

Since $\Pr \{ (U, V) \in \mscr{D} \} = \f{1}{2 \pi} \int \int_{ (u,
v) \in \mscr{D} } \exp \li ( - \f{u^2 + v^2} {2}\ri ) du dv$, using
polar coordinates, we have {\small \bee 2 \pi \Pr \{ (r, \phi) \in
\mscr{D} \} & = & \int \int_{ (r, \phi) \in \mscr{D} } \exp \li ( -
\f{r^2} {2}\ri ) r dr d \phi\\
 & = &  \int \int_{ (r, \phi) \in \mscr{D}^\prime } \exp \li ( -
\f{r^2} {2}\ri ) r dr d \phi + \int \int_{ (r, \phi) \in
\mscr{D}^{\prime \prime} } \exp \li ( - \f{r^2} {2}\ri ) r dr d \phi \\
& = & \int_{\mscr{A}_{\mrm{i}}} \li [ \int_{r = \mcal{B}_{\mrm{v}}
(\phi)}^{\mcal{B}_{\mrm{i}} (\phi)} \exp \li ( - \f{r^2} {2}\ri ) r
dr \ri ] d \phi + \int_{\mscr{A}_{\mrm{v}} \setminus
\mscr{A}_{\mrm{i}}} \li [ \int_{r = \mcal{B}_{\mrm{v}} (\phi)}^\iy
\exp \li ( - \f{r^2} {2}\ri
) r dr \ri ] d \phi\\
& = & \int_{\mscr{A}_{\mrm{i}}} \li [  \exp \li ( -
\f{\mcal{B}_{\mrm{v}}^2 (\phi)} {2}\ri ) -  \exp \li ( -
\f{\mcal{B}_{\mrm{i}}^2 (\phi)} {2}\ri ) \ri ] d \phi +
\int_{\mscr{A}_{\mrm{v}} \setminus \mscr{A}_{\mrm{i}}} \exp \li ( -
\f{\mcal{B}_{\mrm{v}}^2 (\phi)} {2}\ri ) d \phi\\
& = & \int_{\mscr{A}_{\mrm{v}}} \exp \li ( - \f{
\mcal{B}_{\mrm{v}}^2 (\phi) } { 2 } \ri ) d \phi -
\int_{\mscr{A}_{\mrm{i}}} \exp \li ( - \f{ \mcal{B}_{\mrm{i}}^2
(\phi) } { 2 } \ri ) d \phi,
 \eee}
from which the theorem immediately follows.

\section{Proof of Theorem 7}

We use a geometrical approach for proving the theorem.  Let the
horizontal axis be the $u$-axis and the vertical axis be the
$v$-axis.  Note that line $u = k v + g$ intercepts line $u = h$ at
point $R = \li ( h, \f{h - g}{k} \ri )$. Line $u = h$ intercepts the
$u$-axis at $P = (h, 0)$. Line $u = k v + g$ intercepts the $u$-axis
at $Q = (g, 0)$. The theorem can be shown by considering $6$ cases :
(i) $h \leq g < 0$; (ii) $h \leq 0 \leq g$; (iii) $0 < h \leq  g$;
(iv) $0 < g < h$; (v) $g \leq 0 \leq h$; (vi) $ g < h < 0$.

\bsk

In the case of $h \leq g < 0$, $R$ is below the $u$-axis, $P$ is on
the left side of $Q$, and $O$ is on the right side of $Q$. As can be
seen from Figure \ref{fig_QO_Cone_Higher}, the visible and invisible
parts of the  boundary can be expressed, respectively, as {\small
$\mscr{B}_{\mrm{v}} = \li \{ \li ( \f{g}{ \sq{1 + k^2} \cos (\phi +
\phi_k)}, \phi \ri ): \f{\pi}{2} - \phi_k < \phi \leq \pi + \phi_R
\ri \}$} and {\small $\mscr{B}_{\mrm{i}} = \li \{ \li ( \f{h}{\cos
\phi}, \phi \ri ): \f{\pi}{2} < \phi < \pi + \phi_R \ri \}$}. By
Theorem 6 and making use of a change of variable in the integration,
we have $\Pr \{ h \leq U \leq k V + g \} = \int^{\pi + \phi_k +
\phi_R }_{\pi \sh 2} \Psi_{g,k} (\phi) \; d \phi - \int^{\pi +
\phi_R }_{\pi \sh 2} \Psi_h (\phi) \; d \phi$.

\begin{figure}[htbp]
\centerline{\psfig{figure=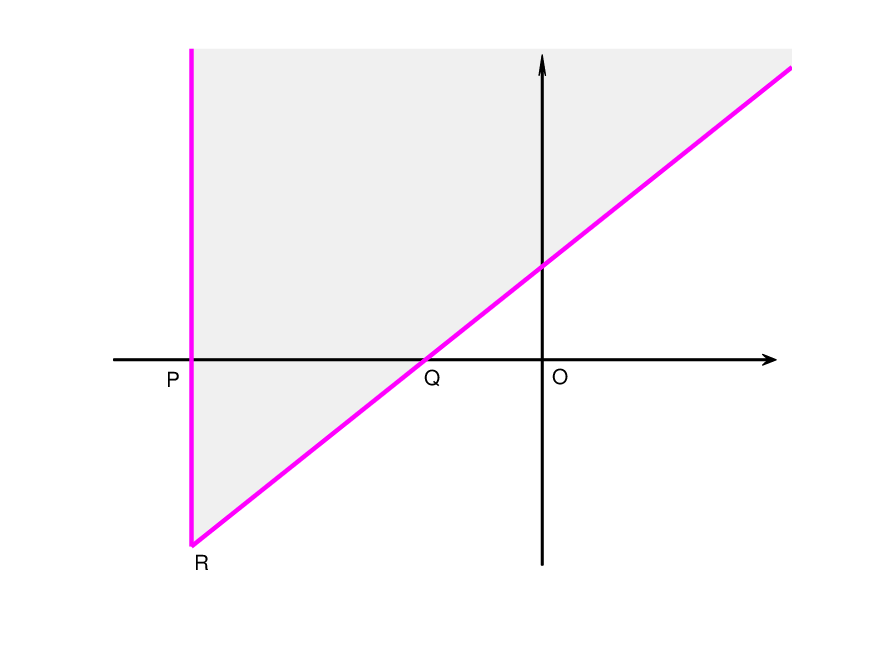, height=2.0in,
width=2.0in }} \caption{Configuration of $h \leq g < 0$ }
\label{fig_QO_Cone_Higher}
\end{figure}

\bsk

In the case of $h \leq 0 \leq g$, $R$ is below the $u$-axis, $P$ is
on the left side of $Q$,  and $O$ is located in between $P$ and $Q$.
As can be seen from Figure \ref{fig_POQ_Cone_Higher}, the boundary
can be expressed as {\small \[ \mscr{B} = \li \{ \li ( \f{h}{\cos
\phi}, \phi \ri ): \f{\pi}{2} < \phi \leq \pi + \phi_R \ri \}
\bigcup \li \{ \li ( \f{g}{\sq{1 + k^2} \cos (\phi + \phi_k)}, \phi
\ri ): \pi + \phi_R \leq \phi < 2 \pi + \f{\pi}{2} - \phi_k \ri \}.
\] } By Theorem 5 and making use of a change of variable in the
integration, we have $\Pr \{ h \leq U \leq k V + g \} = 1 -
\int^{\pi + \phi_R}_{\pi \sh 2}  \Psi_h (\phi) \; d \phi - \int^{3
\pi \sh 2}_{\phi_k + \phi_R}  \Psi_{g,k} (\phi) \; d \phi$.

\begin{figure}[htbp]
\centerline{\psfig{figure=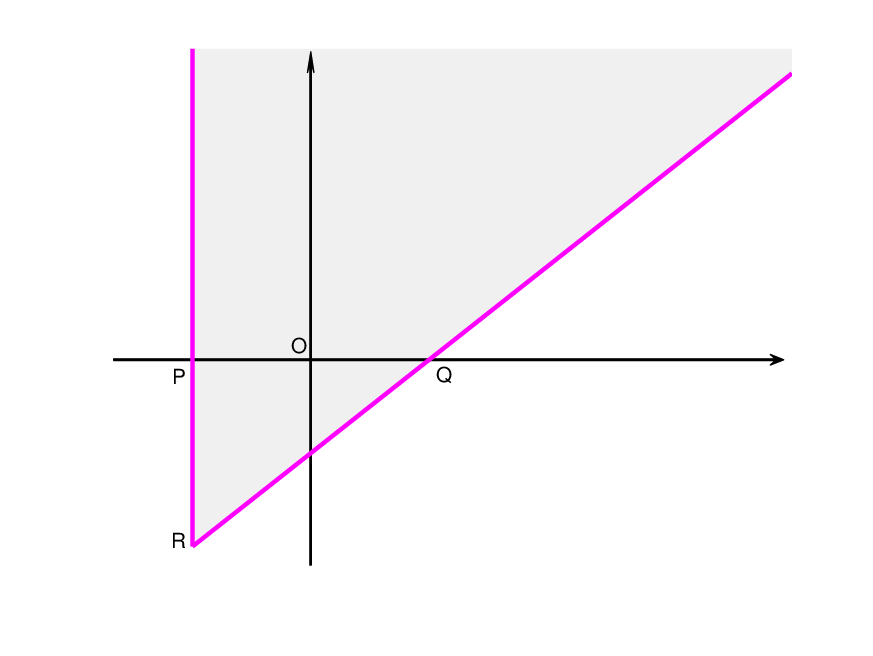, height=2.0in,
width=2.0in }} \caption{Configuration of $h \leq 0 \leq g$ }
\label{fig_POQ_Cone_Higher}
\end{figure}

\bsk

In the case of $0 < h \leq  g$, $O$ is on the left side of $P$, $P$
is on the left side of $Q$, and $R$ is below the $u$-axis.  As can
be seen from Figure \ref{fig_OP_Cone_Higher}, the visible and
invisible parts of the  boundary can be expressed as {\small
$\mscr{B}_{\mrm{v}} = \li \{ \li ( \f{h}{\cos \phi}, \phi \ri ):
\phi_R \leq \phi < \f{\pi}{2} \ri \}$} and {\small
$\mscr{B}_{\mrm{i}} = \li \{ \li ( \f{g}{ \sq{1 + k^2} \cos (\phi +
\phi_k)}, \phi \ri ): \phi_R < \phi < \f{\pi}{2} - \phi_k \ri \}$}
respectively. By Theorem 6 and making use of a change of variable in
the integration, we have
 $\Pr \{ h \leq U \leq k V + g \} = \int^{\pi \sh 2}_{\phi_R} \Psi_h (\phi) \; d \phi - \int^{\pi \sh
2}_{\phi_k + \phi_R} \Psi_{g,k} (\phi) \; d \phi$.

 \begin{figure}[htbp] \centerline{\psfig{figure=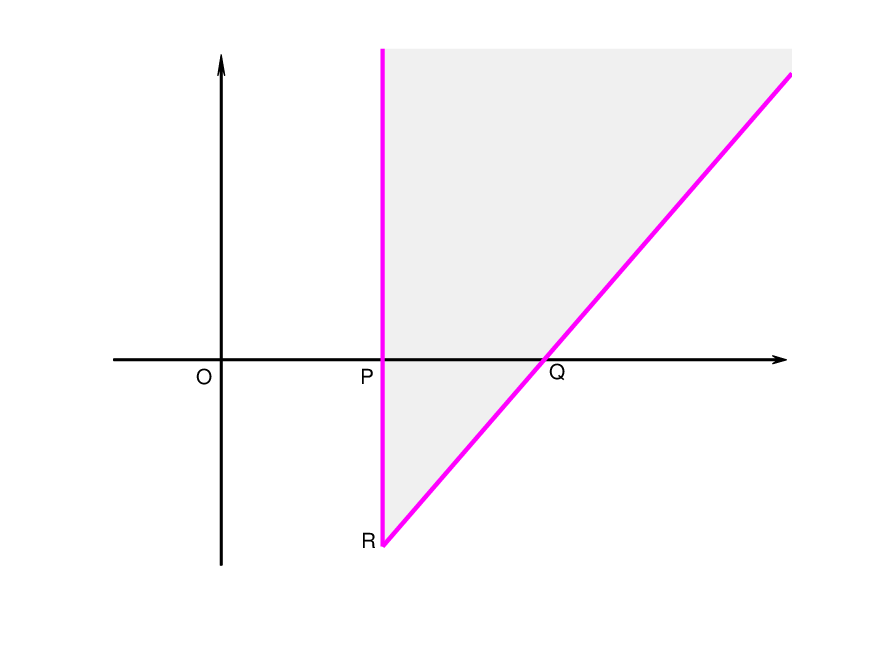,
height=2.0in, width=2.0in }} \caption{Configuration of $0 < h \leq
g$ } \label{fig_OP_Cone_Higher}
\end{figure}

\bsk

In the case of $0 < g < h$, $R$ is above the $u$-axis, $Q$ is on the
left side of $P$, and $O$ is on the left side of $Q$.  As can be
seen from Figure \ref{fig_OQ_Cone_Lower}, the visible and invisible
parts of the  boundary can be expressed as {\small
$\mscr{B}_{\mrm{v}} = \li \{ \li ( \f{h}{\cos \phi}, \phi \ri ):
\phi_R \leq \phi < \f{\pi}{2} \ri \}$} and {\small $
\mscr{B}_{\mrm{i}} = \li \{ \li ( \f{g}{\sq{1 + k^2} \cos (\phi +
\phi_k)}, \phi \ri ): \phi_R < \phi < \f{\pi}{2} - \phi_k \ri \}$}
respectively. By Theorem 6 and making use of a change of variable in
the integration, we have $\Pr \{ h \leq U \leq k V + g \} =
\int^{\pi \sh 2}_{\phi_R} \Psi_h (\phi) \; d \phi - \int^{\pi \sh
2}_{\phi_k + \phi_R} \Psi_{g,k} (\phi) \; d \phi$.

\begin{figure}[htbp]
\centerline{\psfig{figure=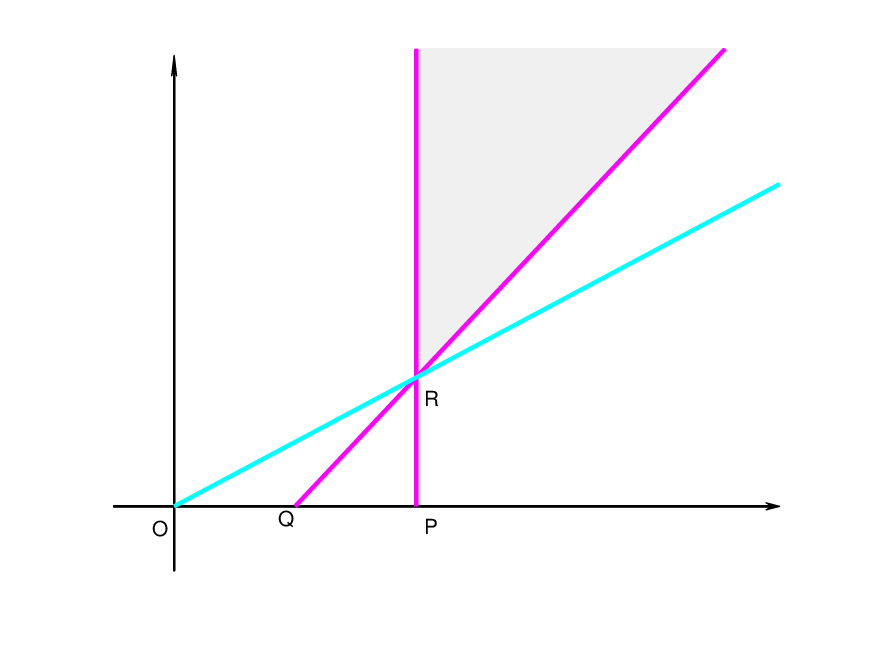, height=2.0in,
width=2.0in }} \caption{Configuration of $0 < g < h$}
\label{fig_OQ_Cone_Lower}
\end{figure}

\bsk

In the case of $g \leq 0 \leq h$, $R$ is above the $u$-axis, $Q$ is
on the left side of $P$, and $O$ is located in between $Q$ and $P$.
As can be seen from Figure \ref{fig_QOP_Cone_Lower}, the boundary is
completely visible and can be expressed as {\small
$\mscr{B}_{\mrm{v}}  = \li \{ \li ( \f{h}{\cos \phi}, \phi \ri ):
\phi_R \leq \phi < \f{\pi}{2} \ri \} \bigcup \li \{ \li (
\f{g}{\sq{1 + k^2} \cos (\phi + \phi_k)}, \phi \ri ): \f{\pi}{2} -
\phi_k < \phi < \phi_R \ri \}$}.  By Theorem 6 and making use of a
change of variable in the integration, we have $\Pr \{ h \leq U \leq
k V + g \} = \int^{\pi \sh 2}_{\phi_R} \Psi_h (\phi) \; d \phi -
\int^{\pi \sh 2}_{\phi_k + \phi_R} \Psi_{g,k} (\phi) \; d \phi$.

\begin{figure}[htbp]
\centerline{\psfig{figure=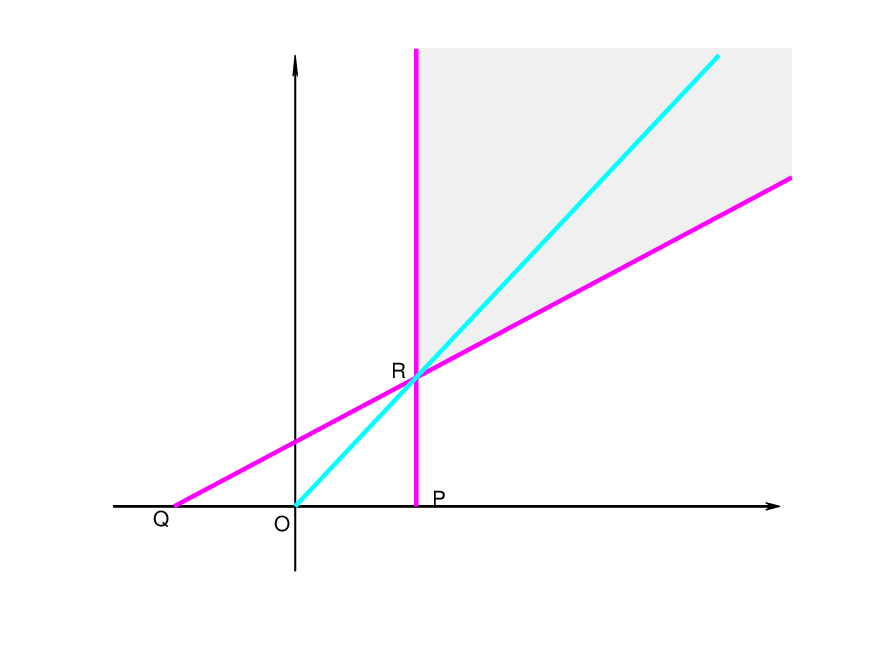, height=2.0in,
width=2.0in }} \caption{Configuration of $g \leq 0 \leq h$}
\label{fig_QOP_Cone_Lower}
\end{figure}

\bsk

In the case of $g < h < 0$, $R$ is above the $u$-axis, $Q$ is on the
left side of $P$, and $P$ is on the left side of $O$. As can be seen
from Figure \ref{fig_PO_Cone_Lower}, the visible and invisible parts
of the boundary can be expressed, respectively, as {\small
$\mscr{B}_{\mrm{v}} = \li \{ \li ( \f{g}{\sq{1 + k^2} \cos (\phi +
\phi_k)}, \phi \ri ): \f{\pi}{2} - \phi_k < \phi \leq \pi + \phi_R
\ri \}$} and {\small $\mscr{B}_{\mrm{i}} = \li \{ \li ( \f{h}{\cos
\phi}, \phi \ri ): \f{\pi}{2} < \phi < \pi + \phi_R \ri \}$}. By
Theorem 6 and making use of a change of variable in the integration,
we have $\Pr \{ h \leq U \leq k V + g \} = \int^{\pi + \phi_k +
\phi_R }_{\pi \sh 2} \Psi_{g,k} (\phi) \; d \phi - \int^{\pi +
\phi_R }_{\pi \sh 2} \Psi_h (\phi) \; d \phi$.  This concludes the
proof of the theorem.

\begin{figure}[htbp]
\centerline{\psfig{figure=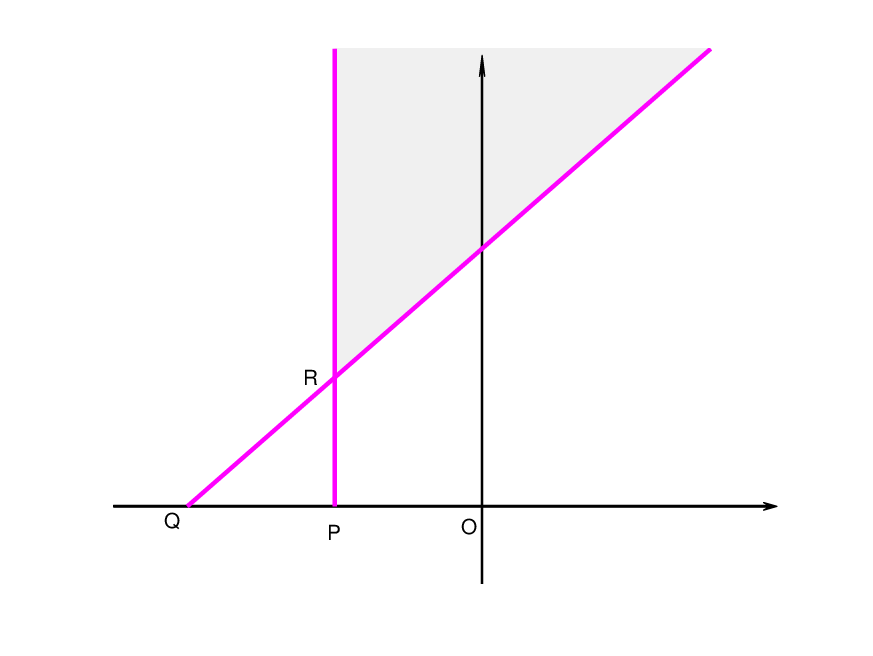, height=2.0in,
width=2.0in }} \caption{Configuration of $g < h < 0$}
\label{fig_PO_Cone_Lower}
\end{figure}

\sect{Proof of Theorem 8}

We shall take a geometrical approach to prove Theorem 8. Before
proceeding to the details of proof, we shall introduce some
notations.  For two points $P_1, \; P_2$ on the $u$-axis, when $P_1$
is on the left side of $P_2$, we write $P_1 < P_2$.  Similarly, when
$P_1$ is on the right side of $P_2$, we write $P_1 >
 P_2$. We use $\wh{P_1 P_2}$ to denote the hyperbolic arc with
end points $P_1$ and $P_2$.  We define some special points $O = (0,
0), \; A = (u_A, v_A), \; B = (u_B, v_B), \; C = (\vse + \sq{h}, 0),
\; D = (\vse - \sq{h}, 0)$ and $M = (\vse, 0)$ that will be
frequently referred in the proof. The domain $\mscr{D}$ is shaded
for all configurations.  The proof of Theorem 8 can be accomplished
by showing Lemmas 2 to 9 in the sequel.

\beL For $\Pr \{ (U, V) \in \mscr{D} \}$ to be non-zero, $\vse, \lm,
g, h, k$ must satisfy one of the following four conditions: (i) $k^2
< \lm, \; g > \sq{h}, \; \vDe \geq 0$; (ii) $k^2 < \lm, \; 0 < g
\leq \sq{h}, \; \vDe \geq 0$; (iii) $k^2 > \lm, \;g k > \sq{ \vDe}$;
(iv) $k^2 > \lm, \;g k \leq \sq{ \vDe}$. \eeL

\bpf Clearly, for $\Pr \{ (U, V) \in \mscr{D} \}$ to be non-zero, a
necessary condition is that there exists at least one tuple $(u, v)$
satisfying equations $\sq{\lm v^2 + h} = u - \vse = k v + g$.  By
letting $z = u - \vse$, we can write the equations as $z - k v = g$
and $(k^2 - \lm) z^2 + 2 \lm g \; z - \lm g^2 - k^2 h = 0$ with $z
\geq 0$, where the discriminant for the quadratic equation of $z$ is
$4 k^2 \vDe$. Therefore, the necessary condition for $\Pr \{ (U, V)
\in \mscr{D} \}$ to be non-zero can be divided as two conditions:
(I) $\vDe \geq 0, \; g \geq 0, \; k^2 < \lm$; (II) $k^2
> \lm$.

If condition (I) holds, then the quadratic equation of $z$ have two
non-negative roots: {\small $z_A = \f{ \lm g - k \sq{ \vDe } } { \lm
- k^2 }, \; z_B = \f{ \lm g + k \sq{ \vDe } } { \lm - k^2 }$}.
Accordingly, there are two tuples $(u_A, v_A)$ and $(u_B, v_B)$
satisfying equations $\sq{\lm v^2 + h} = u - \vse = k v + g$ with
$u_A = z_A + \vse, \; v_A = \f{z_A - g} { k }, \; u_B = z_B + \vse,
\; v_B = \f{z_B - g} { k }$. Noting that $v_A, \; v_B$  are the
roots for equation $(k^2 - \lm ) v^2 + 2 k g v + g^2 - h = 0$ with
respect to $v$, condition (I) can be divided into conditions (i) and
(ii) of the lemma such that (i) implies $\sq{h} + \vse < u_A < u_B,
\; v_A < 0 < v_b$ and that (ii) implies $\sq{h} + \vse < u_A < u_B,
\; 0 \leq v_A < v_B$.

If condition (II) holds, then the quadratic equation of $z$ have two
roots $z_A$ and $z_B$ of opposite signs. Observing that $z_A > z_B$,
we have $z_A = \f{ \lm g - k \sq{ \vDe }   } {\lm - k^2 } > 0 >
z_B$.  Since $v_A = \f{ z_A - g }{ k } = \f{ g k - \sq{ \vDe} } {
\lm - k^2 } \geq 0$ if and only if $g k \leq \sq{ \vDe }$, condition
(II) can be divided into conditions (iii) and (iv) of the lemma such
that (iii) implies $\sq{h} + \vse < u_A, \; v_A < 0$ and that (iv)
implies $\sq{h} + \vse < u_A , \; v_A \geq 0$.  This completes the
proof of the lemma.

\epf

 Now we attempt to express the right branch hyperbola,
 $\mscr{H}_R = \{ (u, v): \sq{\lm v^2 + h} \leq u - \vse \}$ in polar coordinates $(r,
 \phi)$, which is related to the Cartesian coordinates by $u = r \cos \phi, \; v = r \sin \phi$.
Note that the polar coordinates, $(r, \phi)$, of any point of
$\mscr{H}_R$ must satisfy the equation $(r \cos \phi - \vse)^2 - \lm
(r \sin \phi)^2 = h$ with respect to $r \geq 0$, which can be
written as $(\cos^2 \phi - \lm \sin^2 \phi) r^2 - 2 \vse \cos \phi
\; r + \eta = 0$ with $\eta = \vse^2 - h$.  For $\phi$ such that $(h
- \lm \eta) \cos^2 \phi + \lm \eta \geq 0$, we have two real roots
\[
r_\diamond (\phi) = \f{ \eta } { \vse \cos \phi + \sq{ (h - \lm
\eta) \cos^2 \phi + \lm \eta }  }, \qu r_\star (\phi) = \f{ \eta  }
{ \vse \cos \phi - \sq{ (h - \lm \eta) \cos^2 \phi + \lm \eta } } =
- r_\diamond (\phi + \pi).
\]
These are possible expressions for the relationship of polar
coordinates $r$ and $\phi$ of the right branch hyperbola
$\mscr{H}_R$.  However, it is not clear which expression should be
taken. The specific expression and the visibility of $\mscr{H}_R$
are to be determined in the sequel.

\beL \la{lem281} If $O \leq M$, then the right hyperbola
$\mscr{H}_R$ is visible and can be expressed as $\mscr{B}_{\mrm{v}}
= \{ ({r_\star}, \phi) : |\phi| < \phi_{\lm}  \}$. \eeL

\bpf

To show the lemma, we first need to show that ${r_\star} > 0 >
{r_\diamond}$ for $1 - \lm \tan^2 \phi > 0$ and $D < O \leq M$. For
$D < O \leq M$, we have $ \vse - \sq{h} < 0 \leq \vse \Rightarrow
\eta = \vse^2 - h < 0$. Thus, ${r_\diamond}  < 0$ as a result of $1
- \lm \tan^2 \phi > 0 \LRA |\phi| < \phi_{\lm}  < \f{\pi}{2}$. On
the other hand, {\small ${r_\star} = \f{- \eta} { -\vse \cos \phi +
\sq{ (h - \lm \eta) \cos^2 \phi + \lm \eta } }$}. Observing that
$(\vse \cos \phi)^2 - \li [ (h - \lm \eta) \cos^2 \phi + \lm \eta
\ri ] = \eta \cos^2 \phi \; ( 1 - \lm \tan^2 \phi) < 0$ as a
consequence of $\eta < 0$ and $1 - \lm \tan^2 \phi > 0$, we have
${r_\star} > 0$.

Next, we need to show that ${r_\star} > {r_\diamond} \geq 0$ for $1
- \lm \tan^2 \phi > 0$ and $O \leq D$. For $O \leq D$, we have $
\vse - \sq{h} \geq 0 \Rightarrow \eta = \vse^2 - h \geq  0$. Thus,
${r_\diamond}  \geq 0$.  On the other hand, observing that $(\vse
\cos \phi)^2 - \li [ (h - \lm \eta) \cos^2 \phi + \lm \eta \ri ] =
\eta \cos^2 \phi \; ( 1 - \lm \tan^2 \phi) \geq 0$ as a consequence
of $\eta \geq 0$ and $1 - \lm \tan^2 \phi > 0$, we have ${r_\star}
\geq 0$.  Since the denominator of $r_\star$ is smaller than that of
$r_\diamond$, we have ${r_\star} > {r_\diamond} \geq 0$.  This
completes the proof of the lemma.

\epf

\beL \la{lem282} If $M <  O \leq C$, then $\mscr{B}_{\mrm{v}} = \{
({r_\star}, \phi) : |\phi| \leq \phi_{\mrm{m}} \}$ and
$\mscr{B}_{\mrm{i}} = \{ ({r_\diamond}, \phi) : \phi_{\lm}  < |\phi|
< \phi_{\mrm{m}} \}$. \eeL

\bpf Since $M < O \leq C$, we have $ \vse < 0 \leq \vse + \sq{h}
\Rightarrow \eta = \vse^2 - h \leq  0$.  Hence, $(h - \lm \eta)
\cos^2 \phi + \lm \eta  = h \cos^2 \phi + \lm \eta \sin^2 \phi = -
\lm \eta \cos^2 \phi \li ( - \f{h} {\lm \eta} -  \tan^2 \phi \ri )$,
which implies that $(h - \lm \eta) \cos^2 \phi + \lm \eta$ is
nonnegative for $|\phi| \leq \phi_{\mrm{m}}$ and negative for $
\phi_{\mrm{m}} < | \phi | < \f { \pi } {2}$.

 To show the lemma, we first need to show that ${r_\star} \geq 0 \geq {r_\diamond}$  if
$1 - \lm \tan^2 \phi > 0$.  Since $\eta  \leq  0$ and $\vse < 0$, we
have {\small $ {r_\star} = \f{- \eta} { - \vse \cos \phi + \sq{ (h -
\lm \eta) \cos^2 \phi + \lm \eta } } \geq 0$} in view of $1 - \lm
\tan^2 \phi > 0 \LRA |\phi| < \phi_{\lm}  < \f{\pi}{2}$. On the
other hand, observing that {\small ${r_\diamond} = \f{- \eta} {
-\vse \cos \phi - \sq{ (h - \lm \eta) \cos^2 \phi + \lm \eta } }$}
and $(\vse \cos \phi)^2 - \li [ (h - \lm \eta) \cos^2 \phi + \lm
\eta \ri ] = \eta \cos^2 \phi \; ( 1 - \lm \tan^2 \phi) < 0$ as a
consequence of $\eta \leq 0$ and $1 - \lm \tan^2 \phi > 0$, we have
${r_\diamond} \leq 0$.

Next, we need to show that $0 \leq {r_\star} \leq {r_\diamond}$ if
$\phi_{\lm} < |\phi| < \phi_{\mrm{m}} $. By the same argument as
above, we have ${r_\star}  \geq  0$ because $|\phi| < \f{\pi}{2}$.
It remains to show ${r_\star} < {r_\diamond}$. Note that $( \vse
\cos \phi )^2 - \li [ (h - \lm \eta) \cos^2 \phi + \lm \eta \ri ] =
\eta \cos^2 \phi ( 1 - \lm \tan^2 \phi)$ is positive as a result of
$\eta \leq 0$ and $\phi_{\lm} < |\phi| < \phi_{\mrm{m}} \Rightarrow
1 - \lm \tan^2 \phi < 0$. Since $\vse \cos \phi < 0$ as a
consequence of $\vse < 0$ and $\phi_{\lm}  < |\phi| < \phi_{\mrm{m}}
$, it follows that $- \vse \cos \phi - \sq{ (h - \lm \eta) \cos^2
\phi + \lm \eta } > 0$ and thus ${r_\diamond} \geq 0$. Since the
numerators of ${r_\star}$ and ${r_\diamond}$ are equal to the same
non-negative number and the denominator of ${r_\diamond}$ is a
positive number smaller than that of ${r_\star}$, we have
${r_\diamond} \geq {r_\star} \geq 0$.  This completes the proof of
the lemma. \epf

As can be seen from the proof of Lemma 4, the boundary is divided
into visible part $\mscr{B}_{\mrm{v}}$ and invisible part
$\mscr{B}_{\mrm{i}}$ by the upper critical point {\small $ \li ( \f{
\eta } { \vse \cos \phi_m }, \phi_m \ri )$} and the lower critical
point {\small $ \li ( \f{ \eta } { \vse \cos \phi_m }, - \phi_m \ri
)$}.  The visible part is on the left side of the critical line,
which is referred to as the vertical line connecting the lower and
upper critical points. The invisible part is on the right side of
the critical line.

\beL \la{lem284} If $O > C$, then the right hyperbola $\mscr{H}_R$
can be represented as $\{ ({r_\diamond}, \phi) : \phi_{\lm}  < \phi
< 2 \pi - \phi_{\lm} \}$.

\eeL

\bpf

To show the lemma, we first need to show that ${r_\star} < 0 <
{r_\diamond}$ for $\phi_{\lm} < \phi < \pi - \phi_{\lm} $ and $\pi +
\phi_{\lm} < \phi < 2 \phi - \phi_{\lm} $.  Since $O > C$, we have
$\vse < - \sq{h}$ and thus $\eta = \vse^2 - h > 0$.  Since $1 - \lm
\tan^2 \phi < 0$ for $\phi_{\lm}  < \phi < \pi - \phi_{\lm} $ and
$\pi + \phi_{\lm} < \phi < 2 \pi - \phi_{\lm} $, we have $|\vse \cos
\phi| - \sq{ (h - \lm \eta) \cos^2 \phi + \lm \eta } < 0$, leading
to ${r_\star} < 0$. On the other hand, $\vse \cos \phi + \sq{ (h -
\lm \eta) \cos^2 \phi + \lm \eta }
> - |\vse \cos \phi| + \sq{ (h - \lm \eta) \cos^2 \phi + \lm \eta }
> 0$, leading to ${r_\diamond} > 0$.

Next, we need to show that ${r_\star} >  {r_\diamond} > 0$ for $\pi
- \phi_{\lm} < \phi < \pi + \phi_{\lm} $.  For $\pi - \phi_{\lm}  <
\phi < \pi + \phi_{\lm} $, we have $1 - \lm \tan^2 \phi > 0$.  Since
$\eta > 0$ and $\vse < 0$, it must be true that $\vse \cos \phi > 0$
and ${r_\diamond} > 0$.  As a consequence of $\vse \cos \phi > 0$
and $1 - \lm \tan^2 \phi > 0$, we have that the denominator of
${r_\star}$ is positive. Recalling that the numerator of ${r_\star}$
is a positive number $\eta$, we have ${r_\star}
> 0$. Since the numerators of ${r_\star}$ and ${r_\diamond}$ are equal to the same
positive number $\eta$ and the denominator of ${r_\star}$ is  a
positive number smaller than that of ${r_\diamond}$, we have
${r_\star} > {r_\diamond} > 0$. This completes the proof of the
lemma.

 \epf

\beL If $k^2 < \lm, \; g > \sq{h}$ and $\vDe \geq 0$, then $\Pr \{
(U, V) \in \mscr{D} \} = I_{\mrm{np}}$. \eeL

\bpf

As consequence of $k^2 < \lm, \; g > \sq{h}$ and $\vDe \geq 0$, we
have $\sq{h} + \vse < u_A < u_B, \; v_A < 0 < v_B$.  The tangent
line at $A$ intercepts the $u$-axis at $P = (u_P, 0)$ with $u_P$
satisfying {\small $\f{\sq{(u_A - \vse)^2 - h}  } {\sq{\lm} \; (u_A
- u_P) } = \f{ u_A - \vse }{\sq{\lm} \sq{(u_A - \vse)^2 - h}}$},
from which we obtain $u_P = \vse + \f{h}{u_A - \vse} > \vse$.
Similarly, the tangent line at $B$ intercepts the $u$-axis at $Q =
(u_Q, 0)$ with $u_Q = \vse + \f{h}{u_B - \vse} < u_P < u_C$.   Line
$AB$ intercepts the $u$-axis at $R = (u_R, 0)$ with $u_R = g +
\vse$. Clearly, $D < M < Q < P < C$.  The lemma can be shown by
investigating  five cases as follows.

In the case of $\vse + \f{h}{u_B - \vse} \geq 0$, we have $O \leq
Q$. The situation is shown in Figure \ref{fig_OQ_Egg}.  If $O \leq
M$, then, by Lemma 2,  the right branch hyperbola $\mscr{H}_R$ is
completely visible. Accordingly, the visible and invisible parts of
the boundary of $\mscr{D}$ can be expressed, respectively,  as
$\mscr{B}_{\mrm{v}} = \{ ({r_\star}, \phi) : -\phi_A \leq \phi \leq
\phi_B \}$ and $\mscr{B}_{\mrm{i}} = \{ (r_l, \phi) : -\phi_A < \phi
< \phi_B \}$, where {\small $r_l (\phi) = \f{ g + \vse } { \sq{1 +
k^2} \cos (\phi + \phi_k) }$}.  Now consider the situation that $M <
O \leq Q$.  Since the domain, $\mscr{H} = \{ (u, v) : \sq{ \lm v^2 +
h } \leq u - \vse \}$, corresponding to the region included by the
right branch hyperbola $\mscr{H}_R$, is a convex set, we have that
$\mscr{H}$ is divided by line $OA$ into two sub-domains of which one
is below line $OA$ and above the tangent line $PA$, and the other is
above both line $OA$ and the tangent line $PA$.  As can be seen from
Figure \ref{fig_OQ_Egg}, the lower critical point $ \li ( \f{ \eta
}{ \vse \cos \phi_m}, - \phi_m \ri )$ must be below line $OA$. It
follows from Lemma 3 that arc $\wh{AC}$ is visible. By a similar
argument, we have that arc $\wh{CB}$ is visible. Therefore, by Lemma
3, the visible and invisible parts of the boundary of $\mscr{D}$ can
be expressed, respectively, as $\mscr{B}_{\mrm{v}}$ and
$\mscr{B}_{\mrm{i}}$ like the case of $O \leq M$.   Applying Theorem
6 yields {\small $ \Pr \{ (U, V) \in \mscr{D} \} = I_{\mrm{np},1}$}.

\begin{figure}[htbp]
\centerline{\psfig{figure=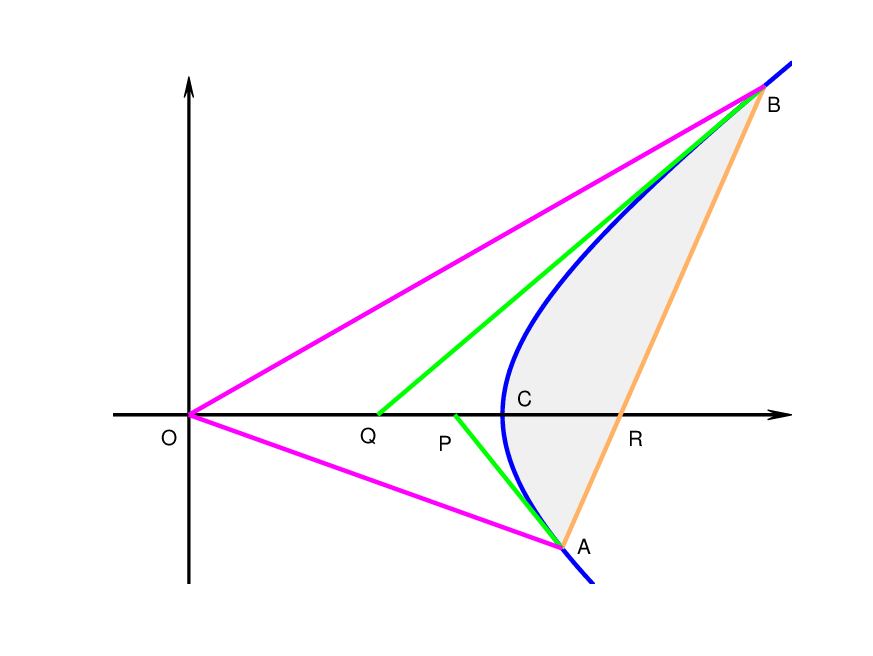, height=1.8in, width=1.8in }}
\caption{Configuration of $O \leq Q$} \label{fig_OQ_Egg}
\end{figure}

In the case of $\vse + \f{h}{u_B - \vse} < 0 \leq \vse + \f{h}{u_A -
\vse}$, we have $Q < O \leq P$. The situation is shown in Figure
\ref{fig_QOP_Egg}. Recall that arc $\wh{AC}$ is visible as in the
preceding case of $O \leq Q$. Since the domain $\mscr{H}$ is a
convex set, we have that $\mscr{H}$ is divided by line $OB$ into two
sub-domains of which one is above line $OB$ and below the tangent
line $QB$, and the other is below both line $OB$ and the tangent
line $QB$.  As can be seen from Figure \ref{fig_QOP_Egg}, the upper
critical point $ \li ( \f{ \eta }{ \vse \cos \phi_m}, \phi_m \ri )$
must be above line $OB$. Hence, applying Lemma 3, the visible and
invisible parts of the boundary of $\mscr{D}$ can be expressed,
respectively,  as {\small $\mscr{B}_{\mrm{v}} = \{ ({r_\star}, \phi)
: -\phi_A \leq \phi \leq \phi_{\mrm{m}} \}$} and {\small
$\mscr{B}_{\mrm{i}} = \{ (r_l, \phi) : -\phi_A < \phi < \phi_B \}
\cup \{ ({r_\diamond}, \phi) : \phi_B \leq \phi < \phi_{\mrm{m}}
\}$}.  Applying Theorem 6 yields {\small $\Pr \{ (U, V) \in \mscr{D}
\} = I_{\mrm{np},2}$}.

\begin{figure}[htbp]
\centerline{\psfig{figure=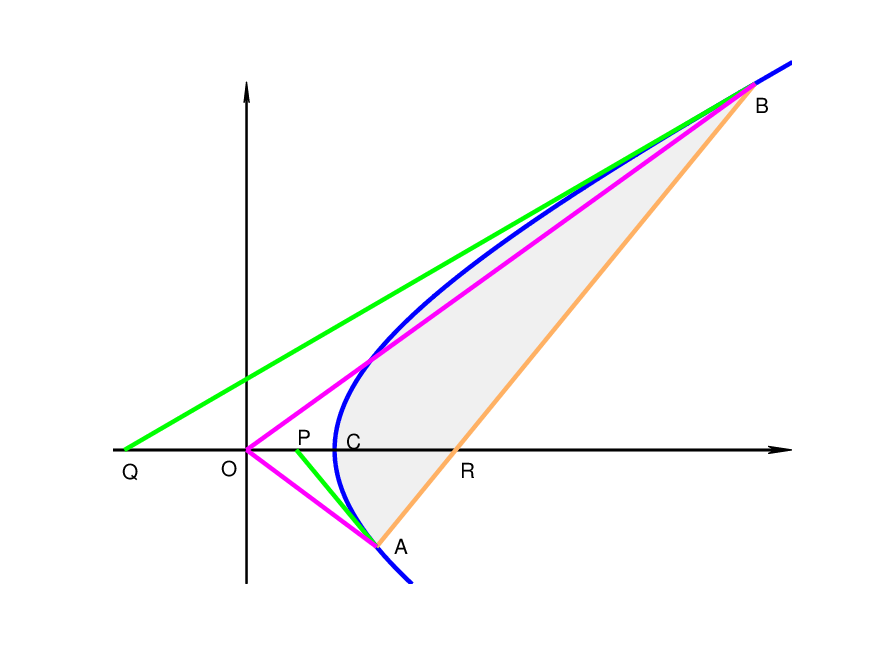, height=1.8in, width=1.8in }}
\caption{Configuration of $Q < O \leq P$} \label{fig_QOP_Egg}
\end{figure}

In the case of $\vse + \f{h}{u_A - \vse} < 0 \leq \vse + \sq{h}$, we
have $P < O \leq C$. The situation is shown in Figure
\ref{fig_POC_Egg}.  By a similar method as that of the case of $Q <
O \leq P$, we have that the upper critical point must be above line
$OB$ and in arc $\wh{CB}$ and that the lower critical point must be
below line $OA$ and in arc $\wh{AC}$.  Hence, by Lemma 3, the
visible and invisible parts of the boundary of $\mscr{D}$ can be
expressed, respectively, as $\mscr{B}_{\mrm{v}} = \{ ({r_\star},
\phi) : -\phi_{\mrm{m}} \leq \phi \leq \phi_{\mrm{m}} \}$ and
$\mscr{B}_{\mrm{i}} = \{ (r_l, \phi) : -\phi_A \leq \phi \leq \phi_B
\} \cup \{ ({r_\diamond}, \phi) : - \phi_{\mrm{m}} < \phi < - \phi_A
\} \cup \{ ({r_\diamond}, \phi) : \phi_B < \phi < \phi_{\mrm{m}}
\}$. By virtue of Theorem 6, we have $\Pr \{ (U, V) \in \mscr{D} \}
= I_{\mrm{np},3}$.

\begin{figure}[htbp]
\centerline{\psfig{figure=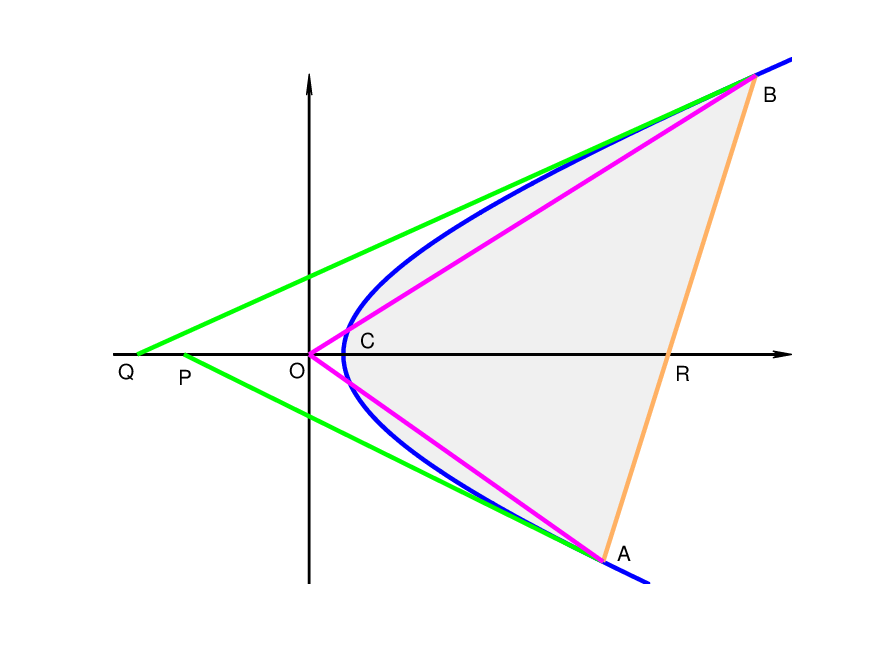, height=1.8in, width=1.8in }}
\caption{Configuration of $P < O \leq C$} \label{fig_POC_Egg}
\end{figure}

In the case of $ \vse + \sq{h} < 0 \leq g + \vse$, we have $C < O
\leq R$. The situation is shown in Figure \ref{fig_COR_Egg}.  By
Lemma 4, the boundary of $\mscr{D}$ can be expressed as $\mscr{B} =
\{ (r_l, \phi) : - \phi_A \leq \phi \leq \phi_B \} \cup \{
({r_\diamond}, \phi) : \phi_B < \phi < 2 \pi - \phi_A \}$. By virtue
of Theorem 5, we have $\Pr \{ (U, V) \in \mscr{D} \} =
I_{\mrm{np},4}$.

\begin{figure}[htbp]
\centerline{\psfig{figure=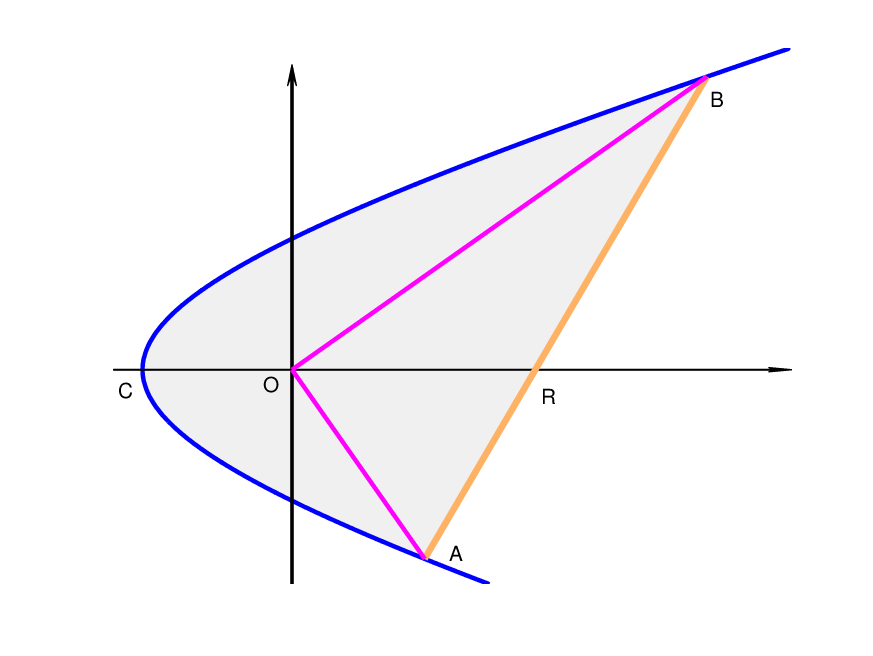, height=1.8in, width=1.8in }}
\caption{Configuration of $C < O \leq R$} \label{fig_COR_Egg}
\end{figure}

In the case of $g + \vse < 0$, we have $O > R$. The situation is
shown in Figure \ref{fig_RO_Egg}. By Lemma 4, the visible and
invisible parts of the boundary of $\mscr{D}$ can be expressed,
respectively, as {\small $\mscr{B}_{\mrm{v}} = \{ (r_l, \phi) :
\phi_B \leq \phi \leq 2 \pi - \phi_A \}$} and {\small
$\mscr{B}_{\mrm{i}} = \{ ({r_\diamond}, \phi) : \phi_B < \phi < 2
\pi - \phi_A \}$}.  By virtue of Theorem 6, we have $\Pr \{ (U, V)
\in \mscr{D} \} = I_{\mrm{np},5}$.

\begin{figure}[htbp]
\centerline{\psfig{figure=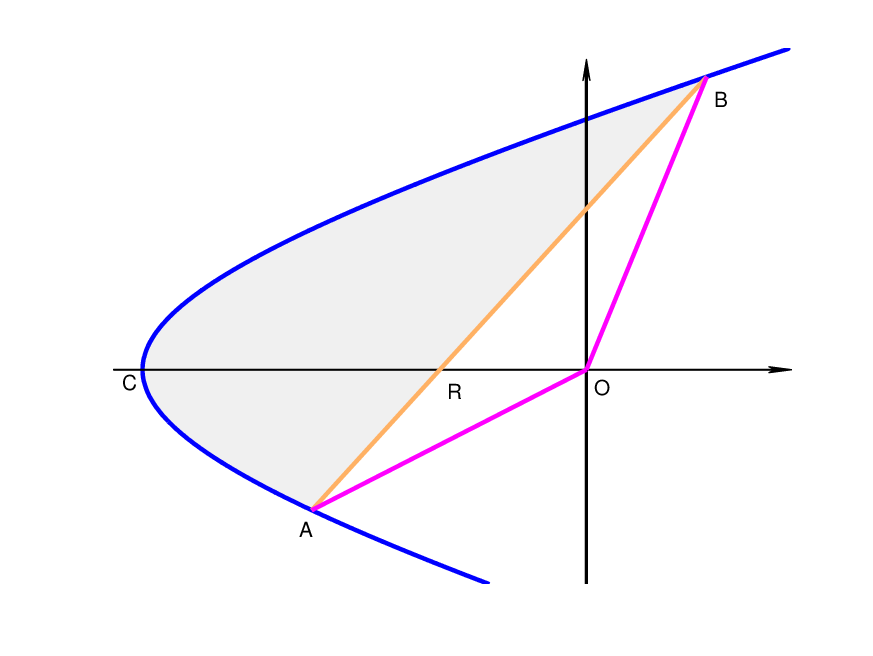, height=1.8in, width=1.8in }}
\caption{Configuration of $O > R$} \label{fig_RO_Egg}
\end{figure}

 \epf

\beL If $k^2 < \lm, \; 0 \leq g \leq \sq{h}$, then $\Pr \{ (U, V)
\in \mscr{D} \} = I_{\mrm{pp}}$.

 \eeL

 \bpf

As a consequence of $k^2 < \lm, \; 0 \leq g \leq \sq{h}$ and $\vDe
\geq 0$, we have $\sq{h} + \vse < u_A < u_B, \; 0 \leq v_A < v_B$.
Clearly, $D < M < Q < R < P < C$. The lemma can be shown by
investigating several cases as follows.

In the case of $\vse + \f{h}{u_B - \vse} \geq 0$, we have $O \leq
Q$. The situation is shown in Figure \ref{fig_OQ_Leave}. By Lemmas 2
and 3,  and a similar argument as that of the first case of Lemma 6,
the visible and invisible parts of the boundary of $\mscr{D}$ can be
determined, respectively, as $\mscr{B}_{\mrm{v}} = \{ ({r_\star},
\phi) : \phi_A \leq \phi \leq \phi_B \}$ and $ \mscr{B}_{\mrm{i}} =
\{ (r_l, \phi) : \phi_A < \phi < \phi_B \}$. By virtue of Theorem 6,
we have $\Pr \{ (U, V) \in \mscr{D} \} = I_{\mrm{pp},1}$.

 \begin{figure}[htbp]
\centerline{\psfig{figure=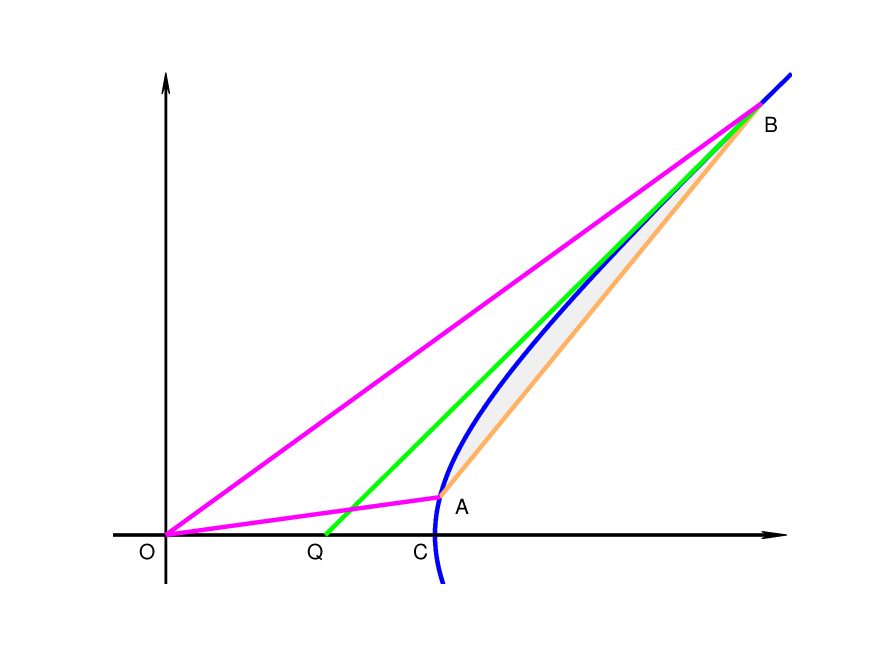, height=1.8in, width=1.8in }}
\caption{Configuration of $O \leq Q$ } \label{fig_OQ_Leave}
\end{figure}

In the case of $\vse + \f{h}{u_B - \vse} < 0 \leq g + \vse$, we have
$Q < O \leq R$. The situation is shown in Figure
\ref{fig_QOR_Leave}.  By Lemma 3 and a similar argument as that of
the second case of Lemma 6, the visible and invisible parts of the
boundary of $\mscr{D}$ can be determined, respectively, as
$\mscr{B}_{\mrm{v}} = \{ ({r_\star}, \phi) : \phi_A \leq \phi \leq
\phi_{\mrm{m}} \}$ and $\mscr{B}_{\mrm{i}} = \{ (r_l, \phi) : \phi_A
< \phi \leq \phi_B \} \cup \{ ({r_\diamond}, \phi) : \phi_B < \phi <
\phi_{\mrm{m}} \}$. By virtue of Theorem 6, we have $\Pr \{ (U, V)
\in \mscr{D} \} = I_{\mrm{pp},2}$.

 \begin{figure}[htbp]
\centerline{\psfig{figure=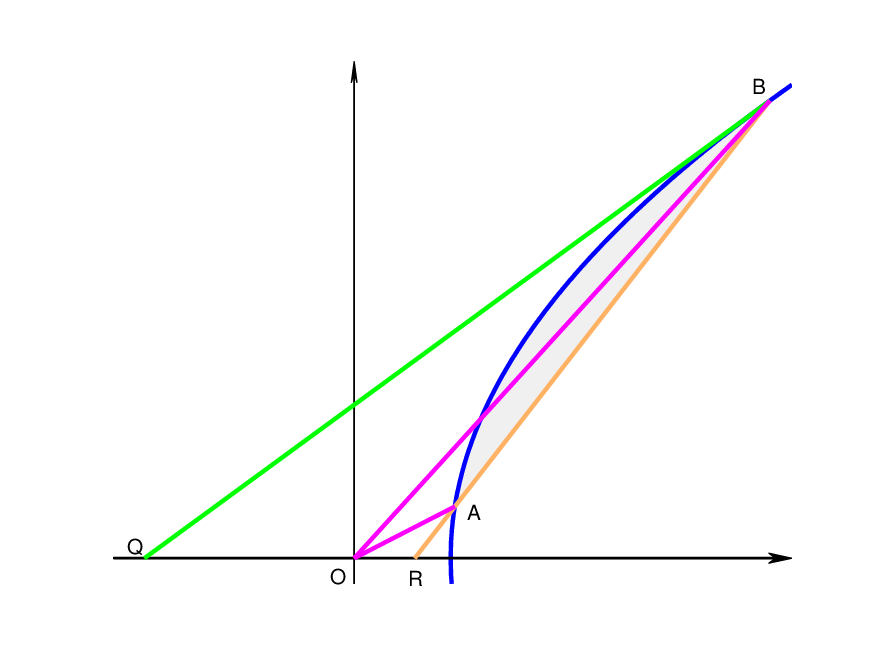, height=1.8in, width=1.8in
}} \caption{Configuration of $Q < O \leq R$} \label{fig_QOR_Leave}
\end{figure}

In the case of $g + \vse < 0 \leq \vse + \f{h}{u_A - \vse}$, we have
$R < O \leq P$. The situation is shown in Figure
\ref{fig_ROP_Leave}.  Observing that the upper critical point must
be above $OA$ and thus must be in arc $\wh{AS}$, by Lemma 3, we have
that the visible and invisible parts of the boundary of $\mscr{D}$
can be expressed, respectively, as $\mscr{B}_{\mrm{v}} = \{
({r_\star}, \phi) : \phi_A \leq \phi \leq \phi_{\mrm{m}} \} \cup \{
(r_l, \phi) : \phi_B \leq \phi < \phi_A \}$ and $\mscr{B}_{\mrm{i}}
= \{ ({r_\diamond}, \phi) : \phi_B < \phi < \phi_{\mrm{m}} \}$. By
virtue of Theorem 6, we have $\Pr \{ (U, V) \in \mscr{D} \} =
I_{\mrm{pp},2}$.

\begin{figure}[htbp]
\centerline{\psfig{figure=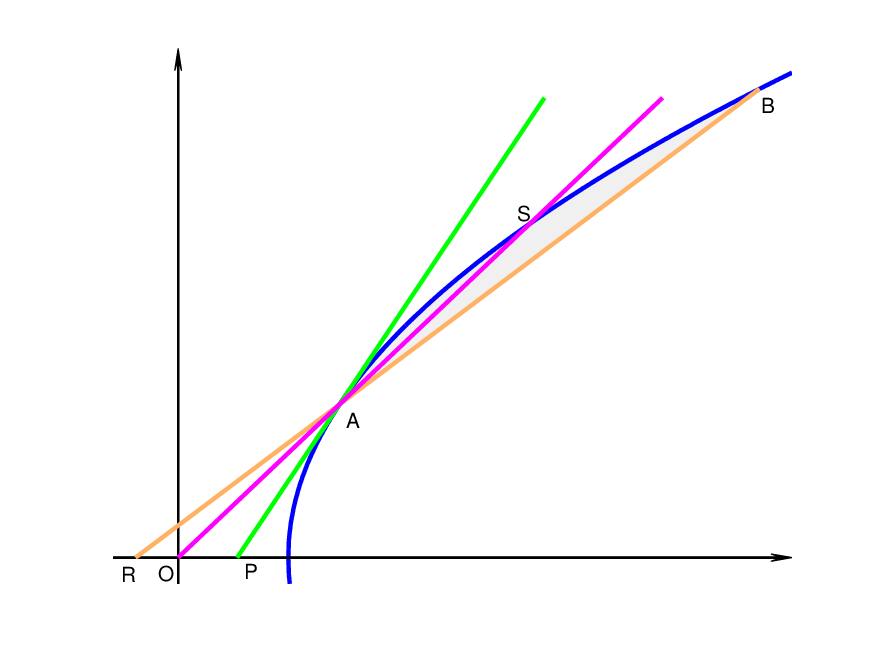, height=2.38in, width=2.38in
}} \caption{Configuration of $R < O \leq P$} \label{fig_ROP_Leave}
\end{figure}

In the case of $\vse + \f{h}{u_A - \vse} < 0 \leq \vse + \sq{h}$, we
have $P < O \leq C$. The situation is shown in Figure
\ref{fig_POC_Leave}.  Observing that the upper critical point must
be in the part of arc $\wh{CA}$ that is above $OA$, by Lemma 3, we
have that the visible and invisible parts of the boundary of
$\mscr{D}$ can be determined, respectively, as $\mscr{B}_{\mrm{v}} =
\{ (r_l, \phi) : \phi_B \leq \phi \leq \phi_A \}$ and
$\mscr{B}_{\mrm{i}} = \{ ({r_\diamond}, \phi) : \phi_B < \phi <
\phi_A \}$. By virtue of Theorem 6, we have $\Pr \{ (U, V) \in
\mscr{D} \} = I_{\mrm{pp},3}$.

\begin{figure}[htbp]
\centerline{\psfig{figure=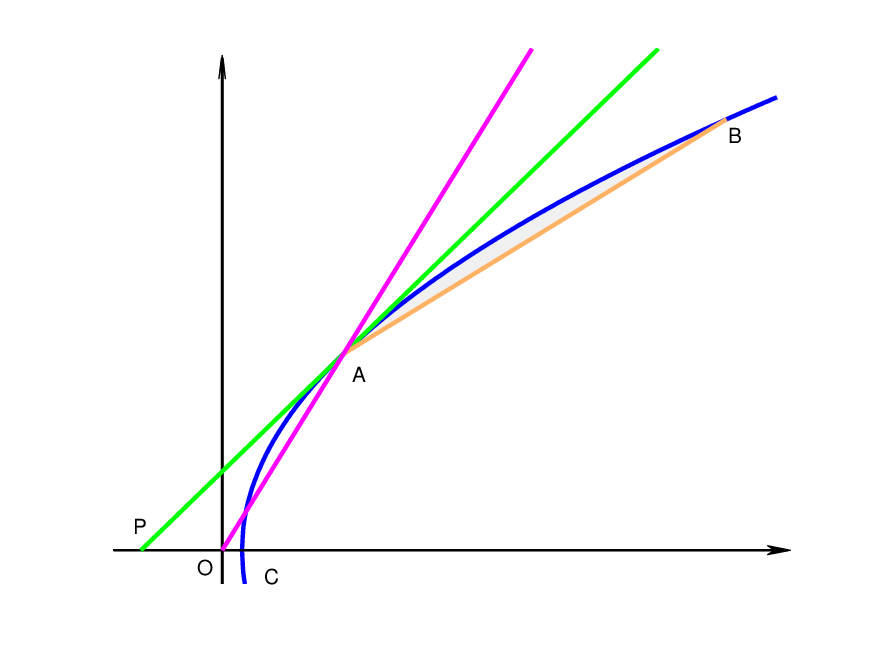, height=1.8in, width=1.8in
}} \caption{Configuration of $P < O \leq C$} \label{fig_POC_Leave}
\end{figure}

In the case of $\vse + \sq{h} < 0$, we have $O > C$. The situation
is shown in Figure \ref{fig_CO_Leave}.  By Lemma 4, the visible and
invisible parts of the boundary of $\mscr{D}$ can be determined,
respectively, as $\mscr{B}_{\mrm{v}} = \{ (r_l, \phi) : \phi_B \leq
\phi \leq \phi_A \}$ and $\mscr{B}_{\mrm{i}} = \{ ({r_\diamond},
\phi) : \phi_B < \phi < \phi_A \}$. By virtue of Theorem 6, we have
$\Pr \{ (U, V) \in \mscr{D} \} = I_{\mrm{pp},3}$.

\begin{figure}[htbp]
\centerline{\psfig{figure=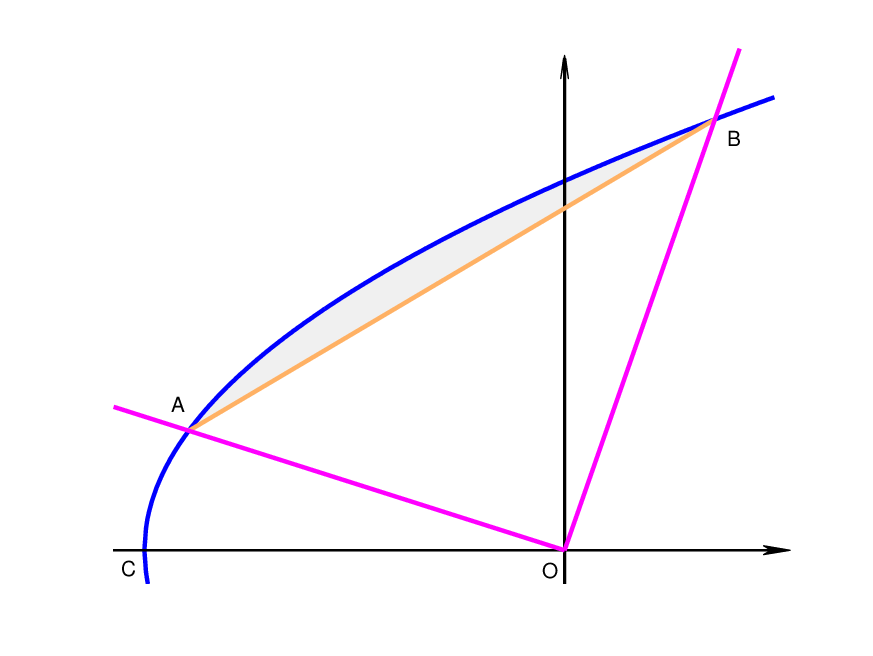, height=1.8in, width=1.8in }}
\caption{Configuration of $O > C$} \label{fig_CO_Leave}
\end{figure}

 \epf

\beL If $k^2 > \lm$ and $g k \leq \sq{ \vDe }$, then $\Pr \{ (U, V)
\in \mscr{D} \} = I_{\mrm{p}}$. \eeL

\bpf

Since $k^2 > \lm$ and $g k \leq \sq{ \vDe }$, we have $v_A \geq 0$.
Consider straight line AB described by equation $u - \vse = k v +
g$, passing through $A = (u_A, v_A)$. Suppose that the tangent line
at $A$ intercepts the $u$-axis at $P$. Draw a line, denoted by $AF$,
from $A$ with angle $\phi_{\lm} $. Extend $FA$ to intercept the
$u$-axis at $G$. Then, $u_A - u_G = \sq{\lm} v_A$, leading to $u_G =
u_A - \sq{\lm} \; v_A$.  The lemma can be shown by considering
several cases as follows.

In the case of $\vse \geq 0$ and $\f{v_A}{u_A} \geq \f{1}{k}$, we
have that  $O \leq M$ and $AB$ is below $OA$. The situation is shown
in Figure \ref{fig_OMG_Tip_1}. Since $O \leq M$, by Lemma 2, the
boundary of $\mscr{D}$ is completely visible and can be expressed as
$\mscr{B}_{\mrm{v}} = \{ ({r_\star}, \phi) : \phi_A \leq \phi <
\phi_{\lm}  \} \cup \li \{ (r_l, \phi) : \f{\pi}{2} - \phi_k < \phi
< \phi_A \ri \}$. By virtue of Theorem 6, we have $\Pr \{ (U, V) \in
\mscr{D} \} = I_{\mrm{p},1}$.

\begin{figure}[htbp]
\centerline{\psfig{figure=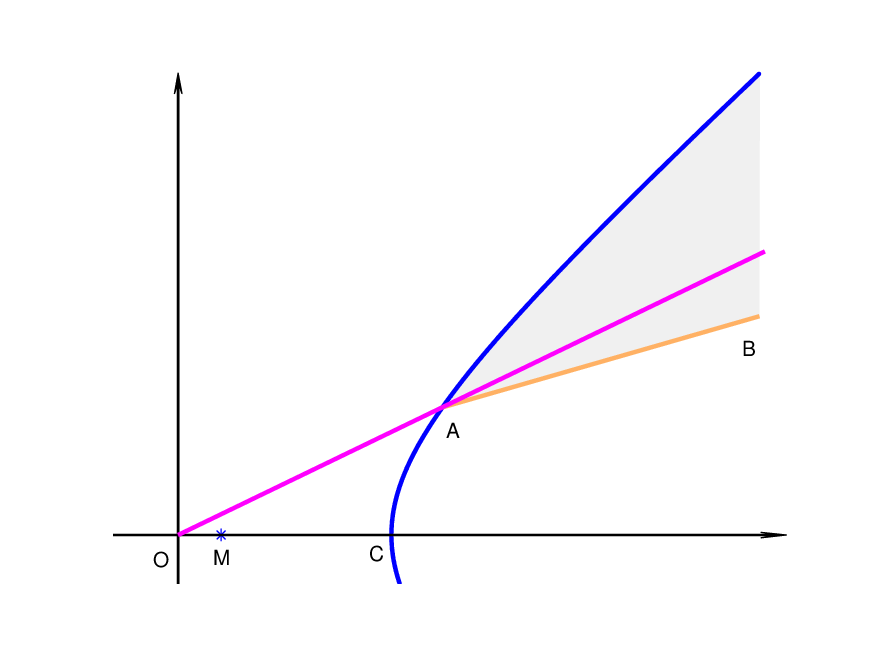, height=1.8in, width=1.8in }}
\caption{Configuration for $O \leq M$ and $AB$ below $OA$ }
\label{fig_OMG_Tip_1}
\end{figure}

In the case of $\vse \geq 0$ and $\f{v_A}{u_A} < \f{1}{k}$, we have
that  $O \leq M$ and $AB$ is above $OA$.  The situation is shown in
Figure \ref{fig_OMG_Tip_2}. By Lemma 2, the visible and invisible
parts of the boundary of $\mscr{D}$ can be determined, respectively,
as $\mscr{B}_{\mrm{v}} = \{ ({r_\star}, \phi) : \phi_A \leq \phi <
\phi_{\lm}  \}$ and $\mscr{B}_{\mrm{i}} = \li \{ (r_l, \phi) :
\phi_A < \phi < \f{\pi}{2} - \phi_k \ri \}$. By virtue of Theorem 6,
we have $\Pr \{ (U, V) \in \mscr{D} \} = I_{\mrm{p},1}$.

\begin{figure}[htbp]
\centerline{\psfig{figure=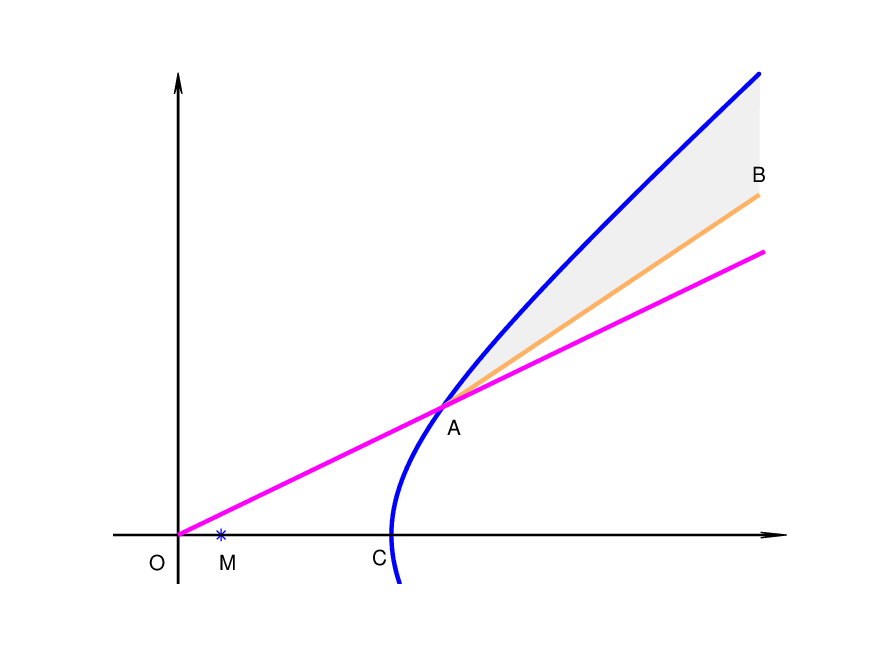, height=1.8in, width=1.8in }}
\caption{Configuration for $O \leq M$ and $AB$ above $OA$ }
\label{fig_OMG_Tip_2}
\end{figure}

In the case of $\vse < 0 \leq u_A - \sq{\lm} \; v_A$ and
$\f{v_A}{u_A} \geq \f{1}{k}$, we have that  $M < O \leq G$ and $AB$
is below $OA$. The situation is shown in Figure \ref{fig_OG_Tip_1}.
Making use of Lemma 3 and the observation that the upper critical
point must be above $OA$, the visible and invisible parts of the
boundary of $\mscr{D}$ can be determined, respectively, as
$\mscr{B}_{\mrm{v}} = \{ ({r_\star}, \phi) : \phi_A \leq \phi \leq
\phi_{\mrm{m}} \} \cup \li \{ (r_l, \phi) : \f{\pi}{2} - \phi_k <
\phi < \phi_A \ri \}$ and $\mscr{B}_{\mrm{i}} = \{ ({r_\diamond},
\phi) : \phi_{\lm} < \phi < \phi_{\mrm{m}} \}$. By virtue of Theorem
6, we have $\Pr \{ (U, V) \in \mscr{D} \} = I_{\mrm{p},2}$.

\begin{figure}[htbp]
\centerline{\psfig{figure=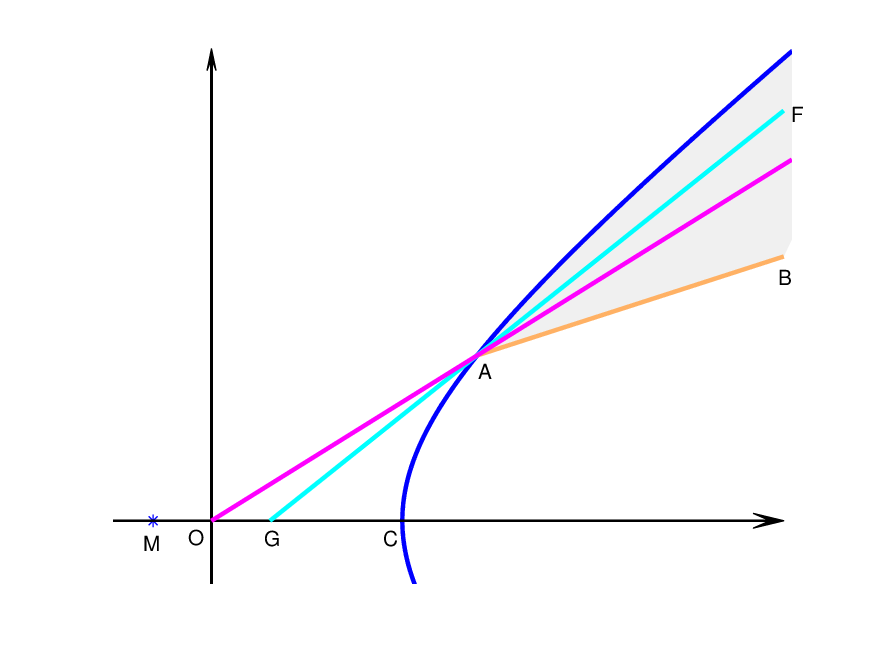, height=1.8in, width=1.8in }}
\caption{Configuration for $M < O \leq G$ and $AB$ below $OA$ }
\label{fig_OG_Tip_1}
\end{figure}

In the case of $\vse < 0 \leq u_A - \sq{\lm} \; v_A$ and
$\f{v_A}{u_A} < \f{1}{k}$, we have that $M < O \leq G$ and $AB$ is
above $OA$.  The situation is shown in Figure \ref{fig_OG_Tip_2}.
Since the upper critical point must be above $OA$, by Lemma 3, the
visible and invisible parts of the boundary of $\mscr{D}$ can be
determined, respectively, as $\mscr{B}_{\mrm{v}} = \{ ({r_\star},
\phi) : \phi_A \leq \phi \leq \phi_{\mrm{m}} \}$ and
$\mscr{B}_{\mrm{i}} = \{ ({r_\diamond}, \phi) : \phi_{\lm}  < \phi <
\phi_{\mrm{m}} \} \cup \li \{ (r_l, \phi) : \phi_A < \phi <
\f{\pi}{2} - \phi_k \ri \}$. By virtue of Theorem 6, we have $\Pr \{
(U, V) \in \mscr{D} \} = I_{\mrm{p},2}$.

\begin{figure}[htbp]
\centerline{\psfig{figure=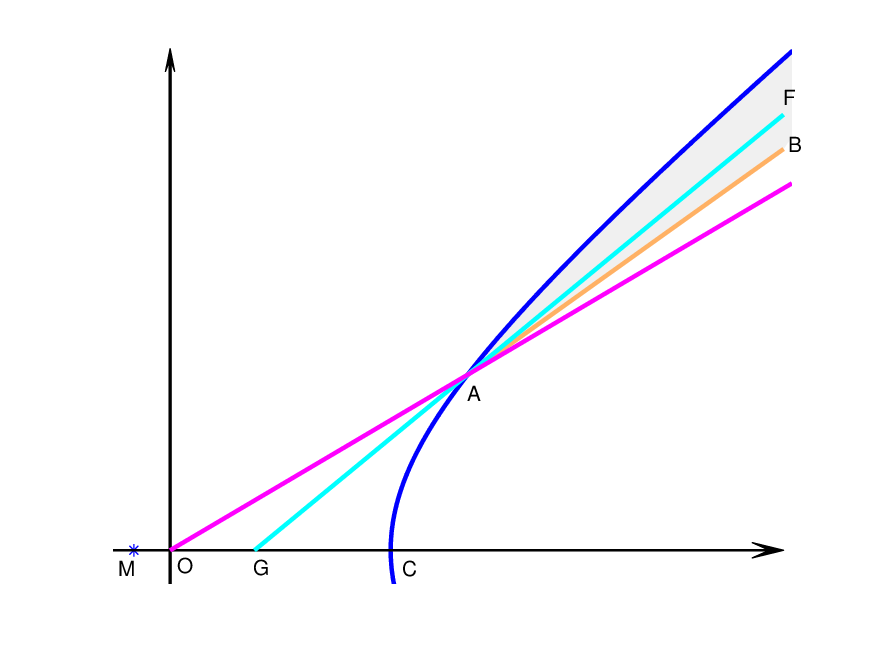, height=1.8in, width=1.8in }}
\caption{Configuration for $M < O \leq G$ and $AB$ above $OA$}
\label{fig_OG_Tip_2}
\end{figure}

In the case of $u_A - \sq{\lm} \; v_A < 0 \leq \vse + \f{h}{u_A -
\vse}$, we have that $G < O \leq P$. The situation is shown in
Figure \ref{fig_GOP_Tip}.   Since $k^2 > \lm$, the slope of line
$AB$ is smaller than that of line $AF$.  As a consequence of $G <
O$, the slope of line $AF$ must be smaller than that of line $OA$.
Hence, the slope of line $AB$ must be smaller than that of line
$OA$. Making use of this observation and noting that the upper
critical point must be above $OA$, we can apply Lemma 3 to determine
the visible and invisible parts of the boundary of $\mscr{D}$,
respectively,  as $\mscr{B}_{\mrm{v}} = \{ ({r_\star}, \phi) :
\phi_A \leq \phi \leq \phi_{\mrm{m}} \} \cup \li \{ (r_l, \phi) :
\f{\pi}{2} - \phi_k < \phi < \phi_A \ri \}$ and $\mscr{B}_{\mrm{i}}
= \{ ({r_\diamond}, \phi) : \phi_{\lm} < \phi < \phi_{\mrm{m}} \}$.
By virtue of Theorem 6, we have $\Pr \{ (U, V) \in \mscr{D} \} =
I_{\mrm{p},2}$.

\begin{figure}[htbp]
\centerline{\psfig{figure=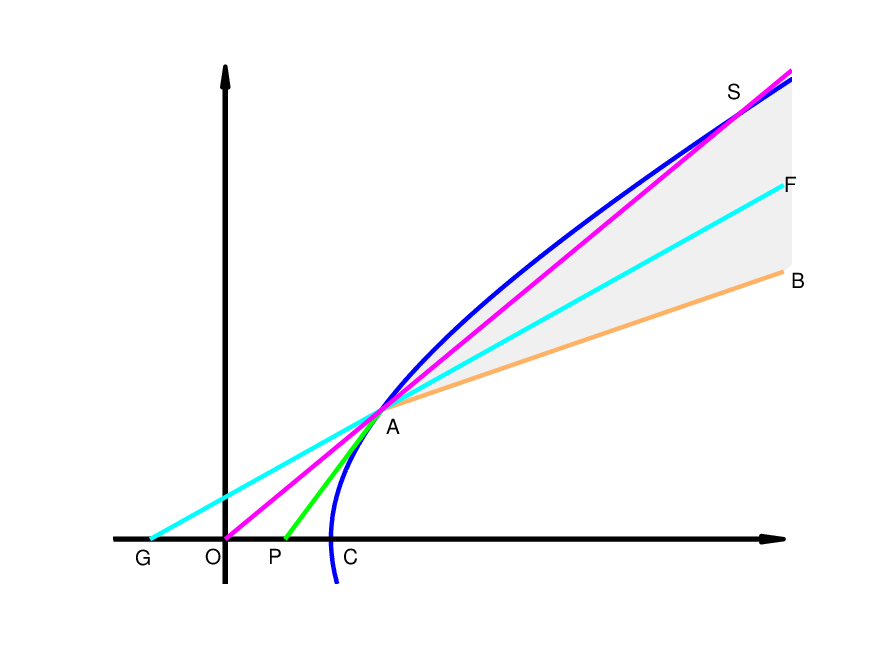, height=1.8in, width=1.8in }}
\caption{Configuration of $G < O \leq P$} \label{fig_GOP_Tip}
\end{figure}

In the case of $\vse + \f{h}{u_A - \vse} < 0 \leq \vse + \sq{h}$, we
have $P < O  \leq C$.  The situation is shown in Figure
\ref{fig_POC_Tip}.  Observing that the upper critical point must be
in the part of arc $\wh{CA}$ that is above line $OA$, by Lemma 3,
the visible and invisible parts of the boundary of $\mscr{D}$ can be
determined, respectively, as $\mscr{B}_{\mrm{v}} = \li \{ (r_l,
\phi) : \f{\pi}{2} - \phi_k < \phi \leq \phi_A \ri \}$ and
$\mscr{B}_{\mrm{i}} = \{ ({r_\diamond}, \phi) : \phi_{\lm} < \phi <
\phi_A \}$. By virtue of Theorem 6, we have $\Pr \{ (U, V) \in
\mscr{D} \} = I_{\mrm{p},3}$.

\begin{figure}[htbp]
\centerline{\psfig{figure=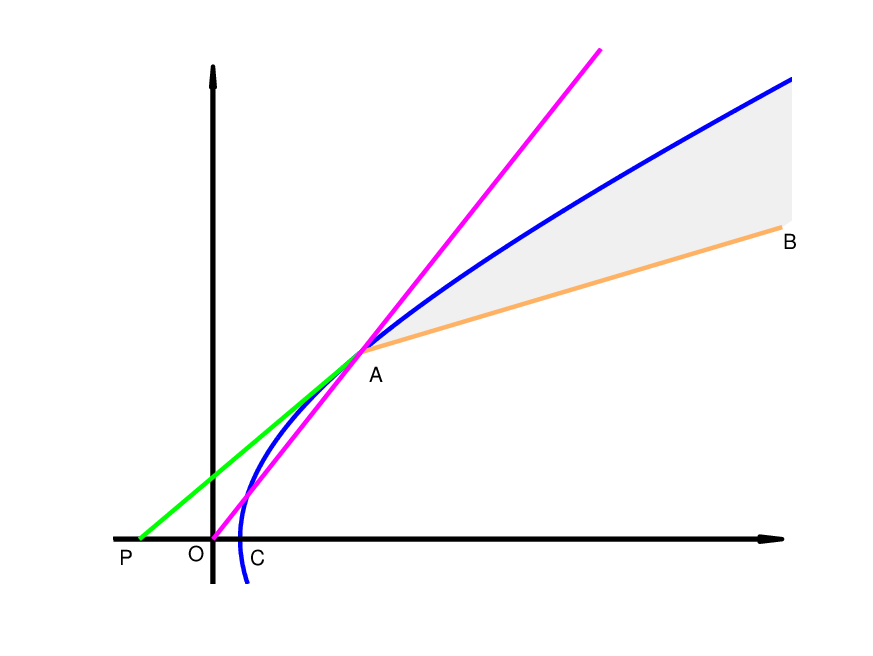, height=1.8in, width=1.8in }}
\caption{Configuration of $P < O \leq C$} \label{fig_POC_Tip}
\end{figure}

In the case of $\vse + \sq{h} < 0$, we have $C < O$.  The situation
is shown in Figure \ref{fig_COR_Tip_Infinite}.  By Lemma 4, the
visible and invisible parts of the boundary of $\mscr{D}$ can be
expressed, respectively, as $\mscr{B}_{\mrm{v}} = \li \{ (r_l, \phi)
: \f{\pi}{2} - \phi_k < \phi \leq \phi_A \ri \}$ and
$\mscr{B}_{\mrm{i}} = \{ ({r_\diamond}, \phi) : \phi_{\lm}  < \phi <
\phi_A \}$. By virtue of Theorem 6, we have $\Pr \{ (U, V) \in
\mscr{D} \} = I_{\mrm{p},3}$.

\begin{figure}[htbp]
\centerline{\psfig{figure=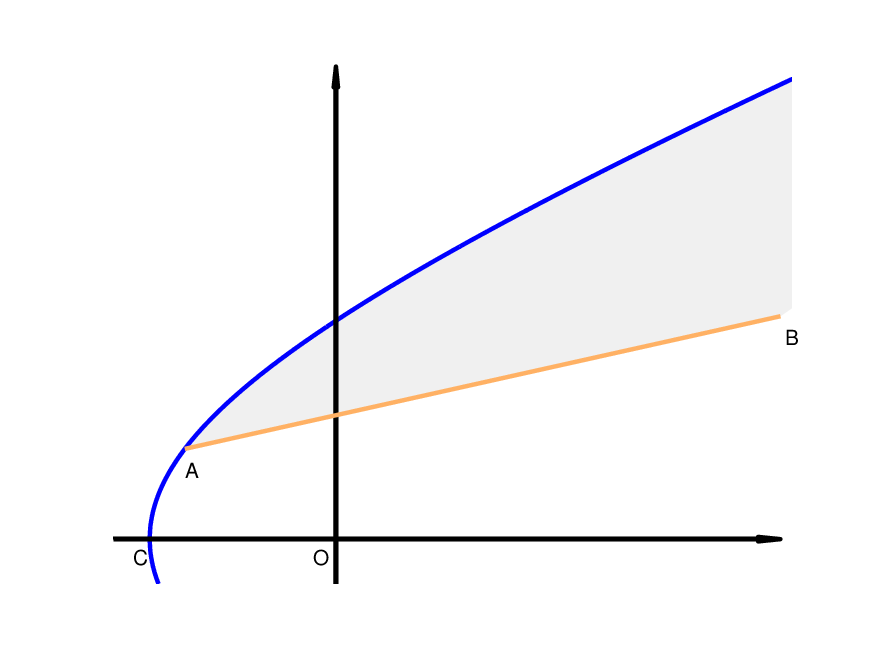, height=1.8in,
width=1.8in }} \caption{Configuration of $C < O$}
\label{fig_COR_Tip_Infinite}
\end{figure}

\epf

\beL

If $k^2 > \lm$ and $g k > \sq{ \vDe}$, then $\Pr \{ (U, V) \in
\mscr{D} \} = I_{\mrm{n}}$.

\eeL

\bpf

For $k^2 > \lm$ and $g k > \sq{ \vDe}$. Then, $v_A < 0$.  The lemma
can be shown by investigating  five cases as follows.

In the case of $\vse \geq 0$, we have $O \leq M$.  The situation is
shown in Figure \ref{OM}.  Since $O \leq M$, by Lemma 2, the right
branch hyperbola $\mscr{H}_R$ is completely visible. Therefore, the
visible and invisible parts of the boundary of $\mscr{D}$ can be
determined, respectively, as $\mscr{B}_{\mrm{v}} = \{ ({r_\star},
\phi) : - \phi_A \leq \phi < \phi_{\lm}  \}$ and $\mscr{B}_{\mrm{i}}
= \li \{ (r_l, \phi) : - \phi_A < \phi < \f{\pi}{2} - \phi_k \ri
\}$. By virtue of Theorem 6, we have $\Pr \{ (U, V) \in \mscr{D} \}
= I_{\mrm{n},1}$.

\begin{figure}[htbp]
\centerline{\psfig{figure=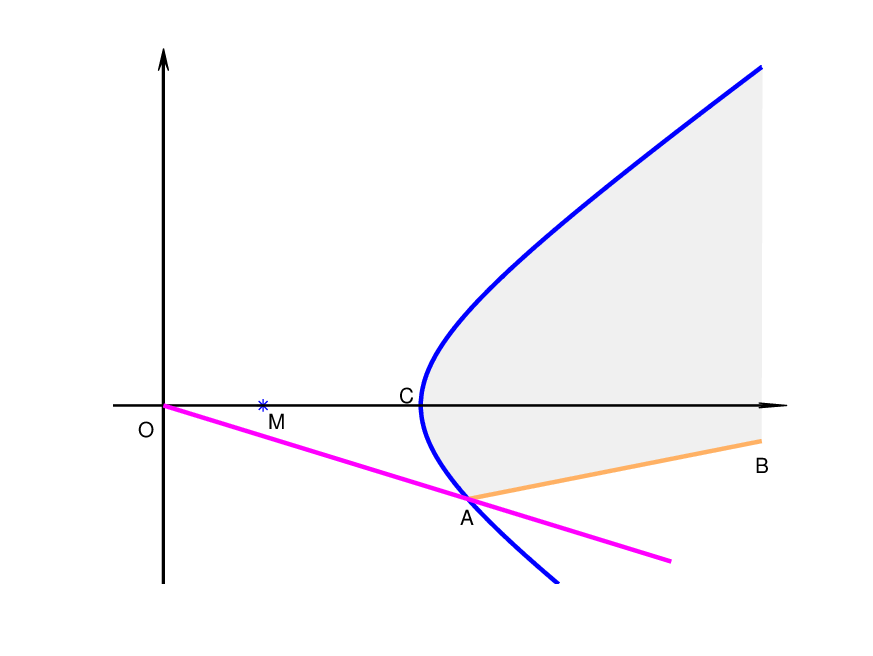, height=1.8in, width=1.8in }}
\caption{Configuration of $O \leq M$} \label{OM}
\end{figure}

In the case of $\vse < 0 \leq \vse + \f{h}{u_A - \vse}$, we have $M
< O \leq P$. The situation is shown in Figure \ref{MOP}. Observing
that the lower critical point must be below line $OA$, by Lemma 3,
we have that arc $\wh{AC}$ must be visible and that the visible and
invisible parts of the boundary of $\mscr{D}$ can be determined,
respectively, as $\mscr{B}_{\mrm{v}} = \{ ({r_\star}, \phi) : -
\phi_A \leq \phi \leq \phi_{\mrm{m}} \}$ and $\mscr{B}_{\mrm{i}} =
\li \{ (r_l, \phi) : - \phi_A < \phi < \f{\pi}{2} - \phi_k \ri \}
\cup \{ ({r_\diamond}, \phi) : \phi_{\lm} < \phi < \phi_{\mrm{m}}
\}$.  By virtue of Theorem 6, we have $\Pr \{ (U, V) \in \mscr{D} \}
= I_{\mrm{n},2}$.

\begin{figure}[htbp]
\centerline{\psfig{figure=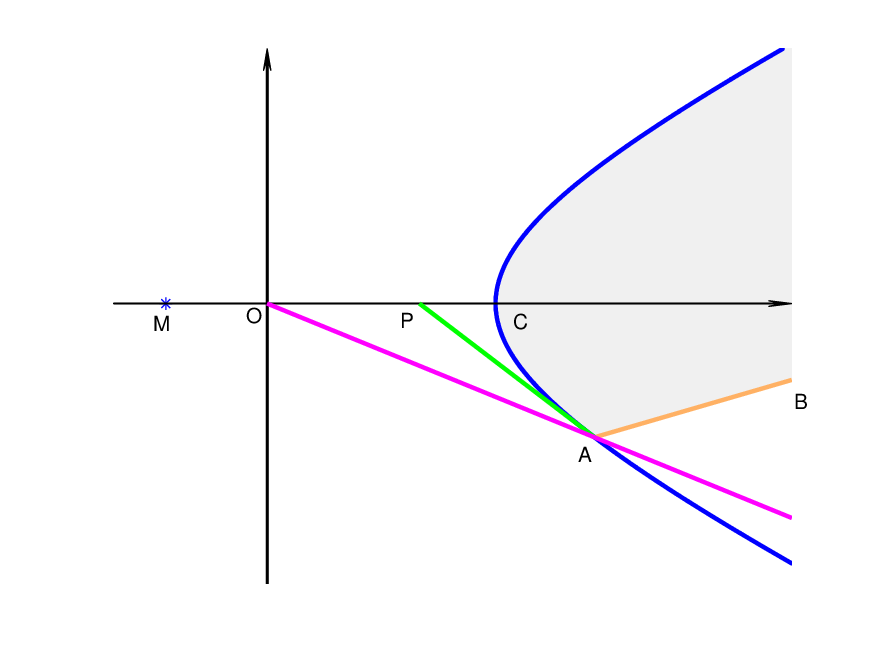, height=1.8in, width=1.8in }}
\caption{Configuration of $M < O \leq P$ } \label{MOP}
\end{figure}

In the case of $\vse + \f{h}{u_A - \vse} < 0 \leq \vse + \sq{h}$, we
have $P < O \leq C$. The situation is shown in Figure \ref{POCm}.
Observing that the lower critical point must be in the part of arc
$\wh{AC}$ that is below line $OA$, by Lemma 3, we have that the
visible and invisible parts of the boundary of $\mscr{D}$ can be
determined, respectively, as $\mscr{B}_{\mrm{v}} = \{ ({r_\star},
\phi) : - \phi_{\mrm{m}} \leq \phi \leq \phi_{\mrm{m}} \}$ and
$\mscr{B}_{\mrm{i}} = \li \{ (r_l, \phi) : - \phi_A < \phi <
\f{\pi}{2} - \phi_k \ri \} \cup \{ ({r_\diamond}, \phi) : \phi_{\lm}
< \phi < \phi_{\mrm{m}} \} \cup \{ ({r_\diamond}, \phi) : -
\phi_{\mrm{m}} < \phi \leq - \phi_A \}$. By virtue of Theorem 6, we
have $\Pr \{ (U, V) \in \mscr{D} \} = I_{\mrm{n},3}$.

\begin{figure}[htbp]
\centerline{\psfig{figure=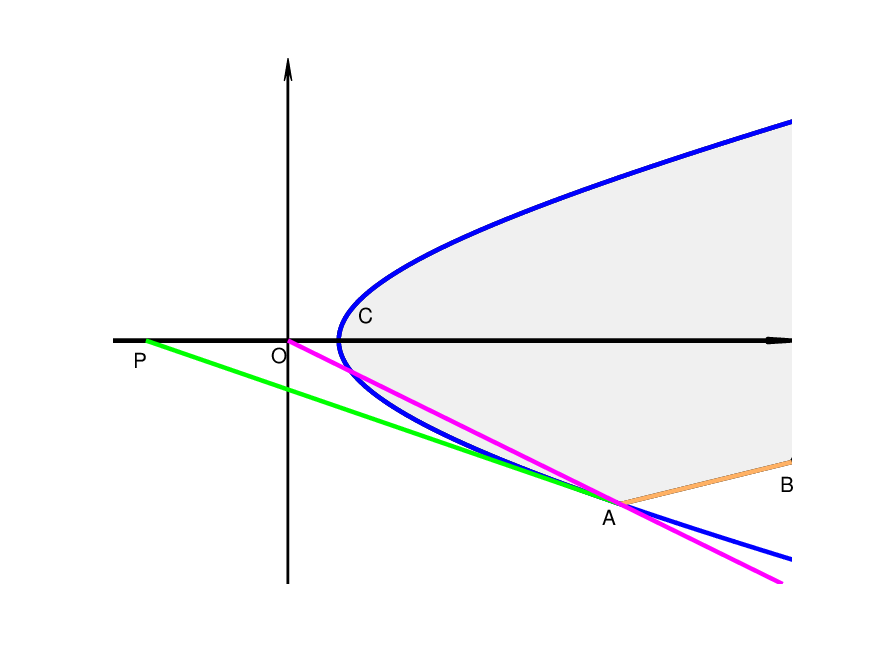, height=1.8in, width=1.8in }}
\caption{Configuration of $P < O \leq C$ } \label{POCm}
\end{figure}

In the case of $\vse + \sq{h} < 0 \leq \vse + g$, we have $C < O
\leq R$.  The situation is shown in Figure \ref{CORm}.  By Lemma 4,
the boundary of $\mscr{D}$ can be expressed as $\mscr{B} = \li \{
(r_l, \phi) : - \phi_A \leq \phi \leq \f{\pi}{2} - \phi_k \ri \}
\cup \{ ({r_\diamond}, \phi) : \phi_{\lm} < \phi < 2 \pi - \phi_A
\}$. By virtue of Theorem 5, we have $\Pr \{ (U, V) \in \mscr{D} \}
= I_{\mrm{n},4}$.

\begin{figure}[htbp]
\centerline{\psfig{figure=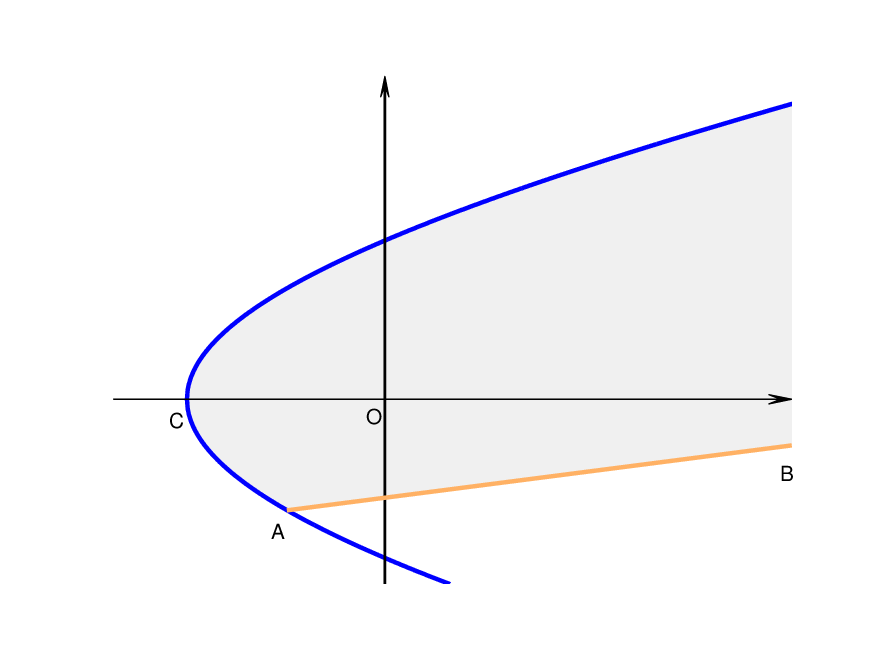, height=1.8in, width=1.8in }}
\caption{Configuration of $C < O \leq R$} \label{CORm}
\end{figure}

In the case of $\vse + g < 0$, we have $R < O$.  The situation is
shown in Figure \ref{ROm}.  By Lemma 4, the visible and invisible
parts of the boundary of $\mscr{D}$ can be determined, respectively,
as $\mscr{B}_{\mrm{v}} = \{ (r_l, \phi) : \f{\pi}{2} - \phi_k < \phi
\leq 2 \pi - \phi_A \}$ and $\mscr{B}_{\mrm{i}} = \{ ({r_\diamond},
\phi) : \phi_{\lm}  < \phi < 2 \pi - \phi_A  \}$. By virtue of
Theorem 6, we have $\Pr \{ (U, V) \in \mscr{D} \} = I_{\mrm{n},5}$.
This completes the proof of the theorem.

\begin{figure}[htbp]
\centerline{\psfig{figure=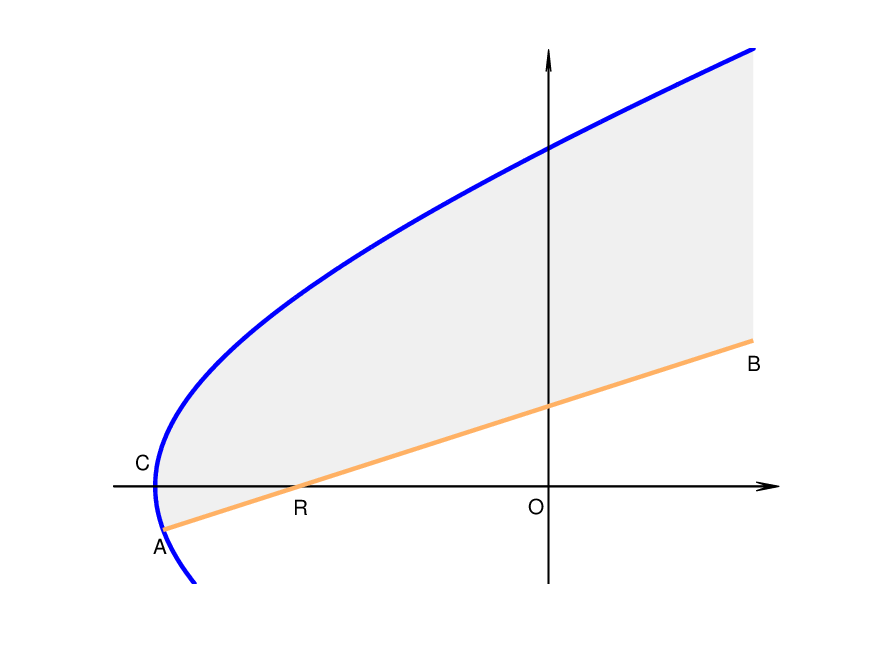, height=1.8in, width=1.8in }}
\caption{Configuration of $R < O$} \label{ROm}
\end{figure}

\epf

\end{document}